\documentclass[11pt]{amsart}
\newtheorem{theorem}{Theorem}[section]
\newtheorem{proposition}[theorem]{Proposition}
\newtheorem{lemma}[theorem]{Lemma}
\newtheorem{corollary}[theorem]{Corollary}
\theoremstyle{definition}
\newtheorem{definition}[theorem]{Definition}

\begin{document}

\title[Ancient solutions to mean curvature flow in higher dimensions]{Uniqueness of convex ancient solutions to mean curvature flow in higher dimensions}
\author{Simon Brendle and Kyeongsu Choi}
\address{Department of Mathematics, Columbia University, 2990 Broadway, New York, NY 10027, USA.}
\address{Department of Mathematics, Massachusetts Institute of Technology, 77 Massachusetts Avenue, Cambridge MA 02138, USA.}
\begin{abstract}
In this paper, we consider noncompact ancient solutions to the mean curvature flow in $\mathbb{R}^{n+1}$ ($n \geq 3$) which are strictly convex, uniformly two-convex, and noncollapsed. We prove that such an ancient solution is a rotationally symmetric translating soliton.
\end{abstract}
\thanks{The first author was supported by the National Science Foundation under grant DMS-1649174 and by the Simons Foundation. The second author was supported by the National Science Foundation under grant DMS-1811267.}
\maketitle

\section{Introduction}

Throughout this paper, we fix an integer $n \geq 3$. Our goal in this paper is to classify all noncompact ancient solutions to mean curvature flow in $\mathbb{R}^{n+1}$ which are convex, uniformly two-convex, and noncollapsed in the sense of Sheng and Wang \cite{Sheng-Wang}: 

\begin{theorem} 
\label{main.thm}
Let $M_t$, $t \in (-\infty,0]$, be a noncompact ancient solution of mean curvature flow in $\mathbb{R}^{n+1}$ which is strictly convex, uniformly two-convex, and noncollapsed. Then $M_t$ is a rotationally symmetric translating soliton.
\end{theorem}

If we evolve a closed, embedded, two-convex hypersurface by mean curvature flow, then it is well known that any blow-up limit is an ancient solution which is weakly convex, uniformly two-convex, and noncollapsed (see \cite{Haslhofer-Kleiner1}, Theorem 1.10, or  \cite{White1},\cite{White2}). If we combine this result with Theorem \ref{main.thm}, we obtain the following result:

\begin{corollary}
Consider an arbitrary closed, embedded, two-convex hypersurface in $\mathbb{R}^{n+1}$, and evolve it by mean curvature flow. At the first singular time, the only possible blow-up limits are shrinking round spheres; shrinking round cylinders; and the unique rotationally symmetric translating soliton.
\end{corollary}

In a recent paper \cite{Brendle-Choi}, we obtained a classification of noncompact ancient solutions in $\mathbb{R}^3$ which are convex and noncollapsed. The proof of Theorem \ref{main.thm} draws on similar techniques. In Section \ref{asymptotic.analysis}, we derive asymptotic estimates for the solution in the cylindrical region. These estimates tell us that, for $-t$ large, the rescaled surface $(-t)^{-\frac{1}{2}} \, M_t \cap B_{5n}(0)$ is $O((-t)^{-\frac{1}{2}})$-close to a cylinder of radius $\sqrt{2(n-1)}$. In Section \ref{barrier}, we combine this estimate with a barrier argument in the spirit of \cite{Angenent-Daskalopoulos-Sesum} to conclude that $\liminf_{t \to -\infty} H_{\text{\rm max}}(t) > 0$, where $H_{\text{\rm max}}(t)$ denotes the supremum of the mean curvature of $M_t$. 

In Section \ref{neck.improvement}, we establish a higher-dimensional version of the Neck Improvement Theorem in \cite{Brendle-Choi}. This step requires significant modifications in the higher-dimensional setting. In order to formulate the Neck Improvement Theorem, we need a notion of $\varepsilon$-symmetry in higher dimensions, which generalizes the one introduced in \cite{Brendle-Choi}. We say that a point $(\bar{x},\bar{t})$ in spacetime is $\varepsilon$-symmetric if there exists a collection of normalized rotation vector fields $\mathcal{K} = \{K_\alpha: 1 \leq \alpha \leq \frac{n(n-1)}{2}\}$ (all having a common axis of rotation) such that $|K_\alpha| \, H \leq 10n$ at $(\bar{x},\bar{t})$ and $|\langle K_\alpha,\nu \rangle| \, H \leq \varepsilon$ in the parabolic neighborhood $\hat{\mathcal{P}}(\bar{x},\bar{t},10,100)$. The main difference between the two-dimensional case and the higher-dimensional case is that, instead of a single rotation vector field in ambient space, we need to consider a collection of normalized rotation vector fields which share a common axis. The statement of the Neck Improvement Theorem can be summarized as follows: if $(\bar{x},\bar{t})$ lies on a neck and every point in a sufficiently large parabolic neighborhood of $(\bar{x},\bar{t})$ is $\varepsilon$-symmetric, then the point $(\bar{x},\bar{t})$ itself is $\frac{\varepsilon}{2}$-symmetric. 

In Section \ref{rotational.symmetry}, we iterate the Neck Improvement Theorem to conclude that any ancient solution which satisfies the assumption of Theorem \ref{main.thm} is rotationally symmetric. Finally, in Section \ref{analysis.in.rotationally.symmetric.case}, we classify ancient solutions with rotational symmetry, thereby completing the proof of Theorem \ref{main.thm}.

\section{Asymptotic analysis as $t \to -\infty$}

\label{asymptotic.analysis}

Suppose that $M_t$, $t \in (-\infty,0]$, is a noncompact ancient solution of mean curvature flow in $\mathbb{R}^{n+1}$ which is strictly convex, uniformly two-convex, and noncollapsed. We consider the rescaled flow $\bar{M}_\tau = e^{\frac{\tau}{2}} \, M_{-e^{-\tau}}$. The hypersurfaces $\bar{M}_\tau$ move with velocity $-(H - \frac{1}{2} \, \langle x,\nu \rangle) \nu$. \\

\begin{proposition}[cf. Haslhofer-Kleiner \cite{Haslhofer-Kleiner1}, Theorem 1.11]
\label{backward.limit}
Consider a sequence $\tau_j \to -\infty$. After passing to a subsequence, the rescaled hypersurfaces $\bar{M}_{\tau_j}$ converge in $C_{\text{\rm loc}}^\infty$ to a cylinder of radius $\sqrt{2(n-1)}$ with axis passing through the origin. 
\end{proposition}

\textbf{Proof.} 
This is proved in \cite{Haslhofer-Kleiner1} on p.~533. We briefly sketch the argument for the convenience of the reader. Let $t_j := -e^{-\tau_j}$. A standard barrier argument implies $\text{\rm dist}(0,M_{t_j}) \leq C(n) \, (-t_j)^{\frac{1}{2}}$ for $j$ sufficiently large. By the speed limit lemma (cf. Lemma 3.4 in \cite{Haslhofer-Kleiner1}), there exist a time $\bar{t}_j \in [t_j,\frac{t_j}{2}]$ and a point $\bar{p}_j \in M_{\bar{t}_j}$ such that $|\bar{p}_j| \leq C(n) \, (-t_j)^{\frac{1}{2}}$ and $H(\bar{p}_j,\bar{t}_j) \leq C(n) \, (-t_j)^{-\frac{1}{2}}$. Using the global curvature estimate in \cite{Haslhofer-Kleiner1} we conclude that, after passing to a subsequence, the rescaled hypersurfaces $\bar{M}_{\tau_j} = (-t_j)^{-\frac{1}{2}} \, M_{t_j}$ converge in $C_{\text{\rm loc}}^\infty$ to a smooth limit. By Huisken's monotonicity formula \cite{Huisken}, the limit must be a self-similar shrinker. Results of Colding and Minicozzi \cite{Colding-Minicozzi} then imply that the limit is either a sphere of radius $\sqrt{2n}$ centered at the origin or a cylinder of radius $\sqrt{2(n-1)}$ with axis passing through the origin. Since the ancient solution $M_t$ is noncompact, the backward limit cannot be a sphere. This completes the proof of Proposition \ref{backward.limit}. \\

In the following, we denote by $\Sigma = \{x \in \mathbb{R}^{n+1}: x_1^2+\hdots+x_n^2=2(n-1)\}$ the cylinder of radius $\sqrt{2(n-1)}$ around the $x_{n+1}$-axis. 

\begin{proposition}
\label{gaussian.area}
For each $\tau$, we have 
\[\int_{\bar{M}_\tau} e^{-\frac{|x|^2}{4}} \leq \int_\Sigma e^{-\frac{|x|^2}{4}}.\]
\end{proposition}

\textbf{Proof.} 
Let us consider an arbitrary sequence $t_j \to -\infty$. The local version of Huisken's monotonicity formula \cite{Huisken} implies that, for each $j$ and each $r>0$, the function 
\[t \mapsto (-t)^{-\frac{n}{2}} \int_{M_t} e^{-\frac{|x|^2}{(-4t)}} \, \Big ( 1 - \frac{|x|^2+2nt}{(-r^2 t_j)} \Big )_+^3\] 
is monotone decreasing for $t \in [t_j,0)$ (see \cite{Ecker}, pp.~64--65). This gives 
\begin{align*} 
&(-t)^{-\frac{n}{2}} \int_{M_t} e^{-\frac{|x|^2}{(-4t)}} \, \Big ( 1 - \frac{|x|^2}{(-r^2 t_j)} \Big )_+^3 \\ 
&\leq (-t_j)^{-\frac{n}{2}} \int_{M_{t_j}} e^{-\frac{|x|^2}{(-4t_j)}} \, \Big ( 1 + \frac{2n}{r^2} - \frac{|x|^2}{(-r^2 t_j)} \Big )_+^3 
\end{align*}
for each $t \in [t_j,0)$ and each $r>0$. We now send $j \to \infty$, keeping $t$ and $r$ fixed. By the monotone convergence theorem, 
\[\lim_{j \to \infty} (-t)^{-\frac{n}{2}} \int_{M_t} e^{-\frac{|x|^2}{(-4t)}} \, \Big ( 1 - \frac{|x|^2}{(-r^2 t_j)} \Big )_+^3 = (-t)^{-\frac{n}{2}} \int_{M_t} e^{-\frac{|x|^2}{(-4t)}}.\] 
Moreover, 
\[\lim_{j \to \infty} (-t_j)^{-\frac{n}{2}} \int_{M_{t_j}} e^{-\frac{|x|^2}{(-4t_j)}} \, \Big ( 1 + \frac{2n}{r^2} - \frac{|x|^2}{(-r^2 t_j)} \Big )_+^3 = \int_\Sigma e^{-\frac{|x|^2}{4}} \, \Big ( 1 + \frac{2n}{r^2} - \frac{|x|^2}{r^2} \Big )_+^3\] 
by Proposition \ref{backward.limit}. Thus, we conclude that 
\[(-t)^{-\frac{n}{2}} \int_{M_t} e^{-\frac{|x|^2}{(-4t)}} \leq \int_\Sigma e^{-\frac{|x|^2}{4}} \, \Big ( 1 + \frac{2n}{r^2} - \frac{|x|^2}{r^2} \Big )_+^3\] 
for each $t$ and each $r >0$. Sending $r \to \infty$, we obtain 
\[(-t)^{-\frac{n}{2}} \int_{M_t} e^{-\frac{|x|^2}{(-4t)}} \leq \int_\Sigma e^{-\frac{|x|^2}{4}}.\] 
This completes the proof of Proposition \ref{gaussian.area}. \\

\begin{lemma}
\label{preliminary.choice.of.rotation.matrix}
We can find a smooth function $\rho(\tau)$ and a function $Q(\tau)$ taking values in $SO(n+1)$ with the following properties: 
\begin{itemize} 
\item $\lim_{\tau \to -\infty} \rho(\tau) = \infty$. 
\item $-\rho(\tau) \leq \rho'(\tau) \leq 0$. 
\item In the ball $B_{2\rho(\tau)}(0)$, the rotated surface $Q(\tau) \bar{M}_\tau$ can be written as a graph over the standard cylinder $\Sigma$, and the $C^4$-norm of the height function is less than $\rho(\tau)^{-8}$.
\end{itemize}
\end{lemma}

\textbf{Proof.} 
We can find a smooth function $Q(\tau)$ taking values in $SO(n+1)$ such that the rotated surfaces $Q(\tau) \bar{M}_\tau$ converge to the standard cylinder $\Sigma$ in $C_{\text{\rm loc}}^\infty$. Hence, we can find a function $\bar{\rho}(\tau)$ with the following properties: 
\begin{itemize} 
\item $\lim_{\tau \to -\infty} \bar{\rho}(\tau) = \infty$. 
\item In the ball $B_{2\bar{\rho}(\tau)}(0)$, the rotated surface $Q(\tau) \bar{M}_\tau$ can be written as a graph over the standard cylinder $\Sigma$, and the $C^4$-norm of the height function is less than $\bar{\rho}(\tau)^{-8}$.
\end{itemize}
We define a function $\rho(\tau)$ by $\rho(\tau) := \inf_{\tau'} \max \{e^{\tau'-\tau},1\} \, \bar{\rho}(\tau')$. Clearly, $\rho(\tau) \leq \bar{\rho}(\tau)$ for each $\tau$, and $\rho(\tau) \to \infty$ as $\tau \to -\infty$. Moreover, the function $\tau \mapsto \rho(\tau)$ is monotone decreasing, and the function $\tau \mapsto e^\tau \, \rho(\tau)$ is monotone increasing. This gives $-\rho(\tau) \leq \rho'(\tau) \leq 0$. Finally, a standard convolution argument allows us to replace the function $\rho(\tau)$ by a smooth function. This completes the proof of Lemma \ref{preliminary.choice.of.rotation.matrix}. \\

As in \cite{Brendle-Choi}, it is necessary to fine tune the choice of the rotation matrix. Let $\varphi \geq 0$ be a smooth cutoff function such that $\varphi=1$ on $[-\frac{1}{2},\frac{1}{2}]$ and $\varphi=0$ outside $[-\frac{2}{3},\frac{2}{3}]$. 

\begin{proposition}
\label{precise.choice.of.rotation.matrix}
Let $\rho(\tau)$ be chosen as in Lemma \ref{preliminary.choice.of.rotation.matrix}. We can find a function $S(\tau)$ taking values in $SO(n+1)$ with the following properties: 
\begin{itemize} 
\item In the ball $B_{2\rho(\tau)}(0)$, the rotated surface $\tilde{M}_\tau := S(\tau) \bar{M}_\tau$ can be written as a graph over the standard cylinder $\Sigma$ of some function $u(\cdot,\tau)$; that is, 
\[\{x + u(x,\tau) \nu_\Sigma(x): x \in \Sigma \cap B_{2\rho(\tau)}(0)\} \subset \tilde{M}_\tau,\] 
where $\nu_\Sigma$ denotes the unit normal to $\Sigma$ and $\|u(\cdot,\tau)\|_{C^4(\Sigma \cap B_{2\rho(\tau)}(0))} \leq \rho(\tau)^{-4}$.
\item The function $u(\cdot,\tau)$ satisfies the orthogonality relations 
\[\int_{\Sigma \cap B_{\rho(\tau)}(0)} e^{-\frac{|x|^2}{4}} \, \langle Ax,\nu_\Sigma \rangle \, u(x,\tau) \, \varphi \Big ( \frac{x_{n+1}}{\rho(\tau)} \Big ) = 0\] 
for every matrix $A \in so(n+1)$. 
\item The matrix $A(\tau) := S'(\tau) S(\tau)^{-1} \in so(n+1)$ satisfies $A(\tau)_{ij} = 0$ for $i,j \in \{1,\hdots,n\}$.
\end{itemize}
\end{proposition}

\textbf{Proof.} 
Given any unit vector $\omega \in S^n$, we consider the rotated cylinder $\Sigma_\omega := \{x \in \mathbb{R}^{n+1}: |x-\langle x,\omega \rangle \omega| = \sqrt{2(n-1)}\}$. Note that $\Sigma_\omega = \Sigma$ if $\omega = (0,\hdots,0,1)$. For each $\omega \in S^n$, we denote by $\pi_\omega: \mathbb{R}^{n+1} \to \Sigma_\omega$ the radial projection to $\Sigma_\omega$. In other words, 
\[\pi_\omega(y) = \langle y,\omega \rangle \omega + \sqrt{2(n-1)} \, \frac{y-\langle y,\omega \rangle \omega}{|y-\langle y,\omega \rangle \omega|}.\] 
For each $\omega \in S^n$, we define a vector $I_\tau(\omega) \in \mathbb{R}^{n+1}$ by 
\begin{align*} 
I_\tau(\omega) 
&:= \int_{\bar{M}_\tau} e^{-\frac{2(n-1)+\langle y,\omega \rangle^2}{4}} \, \Big ( |y-\langle y,\omega \rangle \omega| - \sqrt{2(n-1)} \Big ) \, \varphi \Big ( \frac{\langle y,\omega \rangle}{\rho(\tau)} \Big ) \\ 
&\hspace{20mm} \cdot \det D(\pi_\omega|_{\bar{M}_\tau})(y) \cdot \langle y,\omega \rangle \, \frac{y-\langle y,\omega \rangle \omega}{|y-\langle y,\omega \rangle \omega|}.
\end{align*}
Clearly, $I_\tau(\omega)$ is orthogonal to $\omega$. 

Let $Q(\tau)$ be the rotation matrix defined in Lemma \ref{preliminary.choice.of.rotation.matrix}, and define $\zeta(\tau) \in S^n$ by $Q(\tau) \, \zeta(\tau) = (0,\hdots,0,1)$. Using Lemma \ref{preliminary.choice.of.rotation.matrix}, we obtain 
\[|I_\tau(\zeta(\tau))| \leq C \, \rho(\tau)^{-8}.\] 
Moreover, an analysis of the linearization of the function $I_\tau(\omega)$ near the point $\zeta(\tau)$ gives 
\[|I_\tau(\omega) - I_\tau(\tilde{\omega}) + c(n) \, (\omega-\tilde{\omega})| \leq o(1) \, |\omega-\tilde{\omega}|\] 
if $|\omega - \zeta(\tau)| \leq \rho(\tau)^{-6}$ and $|\tilde{\omega}-\zeta(\tau)| \leq \rho(\tau)^{-6}$. Here, $c(n)$ is a non-zero constant that depends only on $n$. We now apply a standard fixed point theorem to the map which sends $\omega \in S^n$ to $\frac{I_\tau(\omega) + c(n) \, \omega}{|I_\tau(\omega) + c(n) \, \omega|} \in S^n$. Hence, if $-\tau$ is sufficiently large, then we can find a vector $\tilde{\zeta}(\tau) \in S^n$ such that $|\zeta(\tau) - \tilde{\zeta}(\tau)| \leq C \, \rho(\tau)^{-8}$ and 
\[\frac{I_\tau(\tilde{\zeta}(\tau)) + c(n) \, \tilde{\zeta}(\tau)}{|I_\tau(\tilde{\zeta}(\tau)) + c(n) \, \tilde{\zeta}(\tau)|} = \tilde{\zeta}(\tau).\] 
Since $I_\tau(\tilde{\zeta}(\tau))$ is orthogonal to $\tilde{\zeta}(\tau)$, we conclude that $I_\tau(\tilde{\zeta}(\tau)) = 0$. Moreover, $\tilde{\zeta}(\tau)$ depends smoothly on $\tau$.

In the next step, we choose a smooth function $S(\tau)$ taking values in $SO(n+1)$ such that $S(\tau) \, \tilde{\zeta}(\tau) = (0,\hdots,0,1)$. Moreover, we can choose $S(\tau)$ in such a way that the matrix $A(\tau) := S'(\tau) S(\tau)^{-1} \in so(n+1)$ satisfies $A(\tau)_{ij} = 0$ for $i,j \in \{1,\hdots,n\}$. (Otherwise, we replace $S(\tau)$ by $R(\tau)S(\tau)$, where $R(\tau) \in SO(n) \subset SO(n+1)$ is a rotation which fixes the $x_{n+1}$-axis.)

Let $\tilde{M}_\tau := S(\tau) \, \bar{M}_\tau$. Using Lemma \ref{preliminary.choice.of.rotation.matrix} and the estimate $|\zeta(\tau) - \tilde{\zeta}(\tau)| \leq C \, \rho(\tau)^{-8}$, we conclude that the surface $\tilde{M}_\tau$ can be written as a graph over the standard cylinder $\Sigma$ of some function $u(\cdot,\tau)$, where $\|u(\cdot,\tau)\|_{C^4(\Sigma \cap B_{2\rho(\tau)}(0))} \leq \rho(\tau)^{-4}$. Finally, the identity $I_\tau(\tilde{\zeta}(\tau)) = 0$ gives 
\begin{align*} 
&\int_{\tilde{M}_\tau} e^{-\frac{2(n-1)+\langle y,\omega \rangle^2}{4}} \, \Big ( |y-\langle y,\omega \rangle \omega| - \sqrt{2(n-1)} \Big ) \, \varphi \Big ( \frac{\langle y,\omega \rangle}{\rho(\tau)} \Big ) \\ 
&\hspace{20mm} \cdot \det D(\pi_\omega|_{\tilde{M}_\tau})(y) \cdot \langle y,\omega \rangle \, \frac{y-\langle y,\omega \rangle \omega}{|y-\langle y,\omega \rangle \omega|} = 0
\end{align*} 
if $\omega = (0,\hdots,0,1)$ is the vertical unit vector. This finally implies 
\[\int_\Sigma e^{-\frac{2(n-1)+x_{n+1}^2}{4}} \, u(x) \, \varphi \Big ( \frac{x_{n+1}}{\rho(\tau)} \Big ) \, x_{n+1} \, x_i = 0\] 
for $i \in \{1,\hdots,n\}$. From this, the desired orthogonality relations follow. This completes the proof of Proposition \ref{precise.choice.of.rotation.matrix}. \\

We next recall the shrinker foliation constructed in \cite{Angenent-Daskalopoulos-Sesum}. By \cite{Angenent-Daskalopoulos-Sesum}, the union of all the leaves in this foliation contains a truncated cone of the form $\{x \in \mathbb{R}^{n+1}: |x_{n+1}| \geq z_0, \, x_1^2+\hdots+x_n^2 \leq b_0^2 x_{n+1}^2\}$ for some large constant $z_0$ and some small constant $b_0>0$. We denote by $\nu_{\text{\rm fol}}$ the unit normal vector field to this foliation. Moreover, we denote by $\Delta_\tau$ the region in between the cylinder $\Sigma$ and the surface $\tilde{M}_\tau$. 

\begin{proposition}
\label{calibration}
There exists a constant $L_0$ such that for all $L \in [L_0,\rho(\tau)]$ 
\[\int_{\tilde{M}_\tau \cap \{|x_{n+1}| \geq L\}} e^{-\frac{|x|^2}{4}} - \int_{\Sigma \cap \{|x_{n+1}| \geq L\}} e^{-\frac{|x|^2}{4}} \geq -\int_{\Delta_\tau \cap \{|x_{n+1}|=L\}} e^{-\frac{|x|^2}{4}} \, |\langle \omega,\nu_{\text{\rm fol}} \rangle|.\] 
Here, $\omega=(0,\hdots,0,1)$ denotes the vertical unit vector in $\mathbb{R}^{n+1}$.
\end{proposition}

\textbf{Proof.} 
This follows from the identity $\text{\rm div}(e^{-\frac{|x|^2}{4}} \, \nu_{\text{\rm fol}}) = 0$ together with the divergence theorem. See \cite{Brendle-Choi}, Proposition 2.2, for details. \\

\begin{proposition}
\label{consequence.of.gaussian.area.bound}
There exists a constant $L_0$ such that 
\[\int_{\Sigma \cap \{|x_{n+1}| \leq L\}} e^{-\frac{|x|^2}{4}} \, |\nabla^\Sigma u(x,\tau)|^2 \leq C \int_{\Sigma \cap \{|x_{n+1}| \leq \frac{L}{2}\}} e^{-\frac{|x|^2}{4}} \, u(x,\tau)^2\] 
and 
\[\int_{\Sigma \cap \{\frac{L}{2} \leq |x_{n+1}| \leq L\}} e^{-\frac{|x|^2}{4}} \, u(x,\tau)^2 \leq CL^{-2} \int_{\Sigma \cap \{|x_{n+1}| \leq \frac{L}{2}\}} e^{-\frac{|x|^2}{4}} \, u(x,\tau)^2\] 
for all $L \in [L_0,\rho(\tau)]$.
\end{proposition}

\textbf{Proof.} 
Let $\omega=(0,\hdots,0,1)$ denote the vertical unit vector in $\mathbb{R}^{n+1}$. Lemma 4.11 in \cite{Angenent-Daskalopoulos-Sesum} implies that $|\langle \omega,\nu_{\text{\rm fol}} \rangle| \leq C L^{-1} \, |x_1^2+\hdots+x_n^2-2(n-1)|$ for each point $x \in \Delta_\tau \cap \{|x_{n+1}|=L\}$. This gives 
\begin{align*} 
&\int_{\Delta_\tau \cap \{|x_{N+1}|=L\}} e^{-\frac{|x|^2}{4}} \, |\langle \omega,\nu_{\text{\rm fol}} \rangle| \\ 
&\leq C L^{-1} \int_{\Delta_\tau \cap \{|x_{n+1}|=L\}} e^{-\frac{|x|^2}{4}} \, |x_1^2+\hdots+x_n^2-2(n-1)| \\ 
&\leq C L^{-1} \int_{\Sigma \cap \{|x_{n+1}|=L\}} e^{-\frac{|x|^2}{4}} \, u^2. 
\end{align*} 
Combining this estimate with Proposition \ref{calibration} yields 
\[\int_{\tilde{M}_\tau \cap \{|x_{n+1}| \geq L\}} e^{-\frac{|x|^2}{4}} - \int_{\Sigma \cap \{|x_{n+1}| \geq L\}} e^{-\frac{|x|^2}{4}} \geq -C L^{-1} \int_{\Sigma \cap \{|x_{n+1}|=L\}} e^{-\frac{|x|^2}{4}} \, u^2.\] 
We next observe that 
\begin{align*} 
&\int_{\tilde{M}_\tau \cap \{|x_{n+1}| \leq L\}} e^{-\frac{|x|^2}{4}} - \int_{\Sigma \cap \{|x_{n+1}| \leq L\}} e^{-\frac{|x|^2}{4}} \\ 
&= \int_{-L}^L \bigg ( \int_{S^{n-1}} e^{-\frac{z^2}{4}} \, \bigg [ e^{-\frac{(\sqrt{2(n-1)}+u)^2}{4}} \, (\sqrt{2(n-1)}+u)^{n-2} \\ 
&\hspace{40mm} \cdot \sqrt{(\sqrt{2(n-1)}+u)^2 \, \Big ( 1+ \Big ( \frac{\partial u}{\partial z} \Big )^2 \Big ) + |\nabla^{S^{n-1}} u|^2} \\ 
&\hspace{40mm} - e^{-\frac{n-1}{2}} \, \sqrt{2(n-1)}^{n-1} \bigg ] \, \bigg ) \, dz.
\end{align*}
By assumption, the height function $u$ satisfies $|u|+|\frac{\partial u}{\partial z}|+|\nabla^{S^{n-1}} u| \leq o(1)$ for $|x_{n+1}| \leq L$. From this, we deduce that 
\begin{align*} 
&\int_{\tilde{M}_\tau \cap \{|x_{n+1}| \leq L\}} e^{-\frac{|x|^2}{4}} - \int_{\Sigma \cap \{|x_{n+1}| \leq L\}} e^{-\frac{|x|^2}{4}} \\ 
&\geq \int_{-L}^L \bigg ( \int_{S^{n-1}} e^{-\frac{z^2}{4}} \, \Big [ e^{-\frac{(\sqrt{2(n-1)}+u)^2}{4}} \, (\sqrt{2(n-1)}+u)^{n-1} \\ 
&\hspace{40mm} - e^{-\frac{n-1}{2}} \, \sqrt{2(n-1)}^{n-1} + \frac{1}{C} \, |\nabla^\Sigma u|^2 \Big ] \bigg ) \, dz \\ 
&\geq \int_{-L}^L \bigg ( \int_{S^{n-1}} e^{-\frac{z^2}{4}} \, \Big [ -C u^2 + \frac{1}{C} \, |\nabla^\Sigma u|^2 \Big ] \bigg ) \, dz
\end{align*} 
where $C$ is a large constant that depends only on $n$. Putting these facts together, we conclude that
\begin{align*} 
\int_{\tilde{M}_\tau} e^{-\frac{|x|^2}{4}} - \int_\Sigma e^{-\frac{|x|^2}{4}} 
&\geq \int_{\Sigma \cap \{|x_{n+1}| \leq L\}} e^{-\frac{|x|^2}{4}} \, \Big [ -C u^2 + \frac{1}{C} \, |\nabla^\Sigma u|^2 \Big ] \\ 
&- C L^{-1} \int_{\Sigma \cap \{|x_{n+1}|=L\}} e^{-\frac{|x|^2}{4}} \, u^2. 
\end{align*} 
Combining this estimate with Proposition \ref{gaussian.area}, we deduce that 
\begin{align*} 
&\int_{\Sigma \cap \{|x_{n+1}| \leq L\}} e^{-\frac{|x|^2}{4}} \, |\nabla^\Sigma u|^2 \\
&\leq C \int_{\Sigma \cap \{|x_{n+1}| \leq L\}} e^{-\frac{|x|^2}{4}} u^2 + C L^{-1} \int_{\Sigma \cap \{|x_{n+1}|=L\}} e^{-\frac{|x|^2}{4}} \, u^2. 
\end{align*}
On the other hand, using the divergence theorem, we obtain
\begin{align*} 
L \int_{\Sigma \cap \{|x_{n+1}|=L\}} e^{-\frac{|x|^2}{4}} \, u^2 
&= \int_{\Sigma \cap \{|x_{n+1}| \leq L\}} \text{\rm div}_\Sigma(e^{-\frac{|x|^2}{4}} \, u^2 \, x^{\text{\rm tan}}) \\ 
&= \int_{\Sigma \cap \{|x_{n+1}| \leq L\}} e^{-\frac{|x|^2}{4}} \, \Big ( u^2 - \frac{1}{2} \, x_{n+1}^2 \, u^2 + 2u \, \langle x^{\text{\rm tan}},\nabla^\Sigma u \rangle \Big ) \\ 
&\leq \int_{\Sigma \cap \{|x_{n+1}| \leq L\}} e^{-\frac{|x|^2}{4}} \, \Big ( u^2 - \frac{1}{4} \, x_{n+1}^2 \, u^2 + 4 \, |\nabla^\Sigma u|^2 \Big ),
\end{align*}
and consequently 
\begin{align*} 
&L^2 \int_{\Sigma \cap \{|x_{n+1}| \leq L\}} e^{-\frac{|x|^2}{4}} u^2 + L \int_{\Sigma \cap \{|x_{n+1}|=L\}} e^{-\frac{|x|^2}{4}} \, u^2 \\ 
&\leq C \int_{\Sigma \cap \{|x_{n+1}| \leq L\}} e^{-\frac{|x|^2}{4}} \, |\nabla^\Sigma u|^2 + C L^2 \int_{\Sigma \cap \{|x_{n+1}| \leq \frac{L}{2}\}} e^{-\frac{|x|^2}{4}} \, u^2. 
\end{align*}
To summarize, we have shown that 
\begin{align*} 
&\int_{\Sigma \cap \{|x_{n+1}| \leq L\}} e^{-\frac{|x|^2}{4}} \, |\nabla^\Sigma u|^2 \\
&\leq CL^{-2} \int_{\Sigma \cap \{|x_{n+1}| \leq L\}} e^{-\frac{|x|^2}{4}} \, |\nabla^\Sigma u|^2 + C \int_{\Sigma \cap \{|x_{n+1}| \leq \frac{L}{2}\}} e^{-\frac{|x|^2}{4}} \, u^2. 
\end{align*} 
If $L$ is sufficiently large, we can absorb the first term on the right hand side into the left hand side. This gives 
\[\int_{\Sigma \cap \{|x_{n+1}| \leq L\}} e^{-\frac{|x|^2}{4}} \, |\nabla^\Sigma u|^2 \leq C \int_{\Sigma \cap \{|x_{n+1}| \leq \frac{L}{2}\}} e^{-\frac{|x|^2}{4}} \, u^2.\] 
This proves the first statement. Using the inequality 
\[0 \leq \int_{\Sigma \cap \{|x_{n+1}| \leq L\}} e^{-\frac{|x|^2}{4}} \, \Big ( u^2 - \frac{1}{4} \, x_{n+1}^2 \, u^2 + 4 \, |\nabla^\Sigma u|^2 \Big ),\] 
the second statement follows. This completes the proof of Proposition \ref{consequence.of.gaussian.area.bound}. \\

Let us denote by $\mathcal{H}$ the space of all functions $f$ on $\Sigma$ such that 
\[\|f\|_{\mathcal{H}}^2 = \int_\Sigma e^{-\frac{|x|^2}{4}} \, f^2 < \infty.\] 
We define an operator $\mathcal{L}$ on the cylinder $\Sigma$ by 
\[\mathcal{L} f = \Delta_\Sigma f - \frac{1}{2} \, \langle x^{\text{\rm tan}},\nabla^\Sigma f \rangle + f.\] 
This operator can be rewritten as 
\[\mathcal{L} f = \frac{\partial^2}{\partial z^2} f + \frac{1}{2(n-1)} \, \Delta_{S^{n-1}} f - \frac{1}{2} \, z \, \frac{\partial}{\partial z} f + f.\] 
Let $Y_m$ be a basis of eigenfunctions of $\Delta_{S^{n-1}}$, and let $\lambda_m$ denote the corresponding eigenvalues. We assume that the eigenfunctions $Y_m$ are normalized so that $\int_{S^{n-1}} Y_m(\theta)^2 \, d\theta = \frac{1}{n} \, |S^{n-1}|$. Note that $\lambda_0=0$, $\lambda_1=\hdots=\lambda_n=n-1$, and $\lambda_{n+1}=2n$. Moreover, $Y_0(\theta)=\frac{1}{\sqrt{n}}$, $Y_1(\theta)=\theta_1, \hdots, Y_n(\theta)=\theta_n$ for $\theta \in S^{n-1}$, where $\theta_1,\hdots,\theta_n$ denote the Cartesian coordinates of $\theta$.

The eigenfunctions of $\mathcal{L}$ are of the form $H_l \big ( \frac{z}{2} \big ) \, Y_m$, where $H_l$ denotes the Hermite polynomial of degree $l$. The corresponding eigenvalues are given by $1-\frac{l}{2}-\frac{\lambda_m}{2(n-1)}$. Thus, there are $n+2$ eigenfunctions that correspond to positive eigenvalues of $\mathcal{L}$, and these are given by $1,x_{n+1},x_1, \hdots, x_n$, up to scaling. The span of these eigenfunctions will be denoted by $\mathcal{H}_+$. Moreover, there are $n+1$ eigenfunctions of $\mathcal{L}$ with eigenvalue $0$, and these are given by $x_{n+1}^2-2, x_1 x_{n+1}, \hdots,x_n x_{n+1}$, up to scaling. The span of these eigenfunctions will be denoted by $\mathcal{H}_0$. The span of all remaining eigenfunctions will be denoted by $\mathcal{H}_-$. With this understood, we have 
\begin{align*} 
&\langle \mathcal{L} f,f \rangle_{\mathcal{H}} \geq \frac{1}{2} \, \|f\|_{\mathcal{H}}^2 & \text{\rm for $f \in \mathcal{H}_+$,} \\ 
&\langle \mathcal{L} f,f \rangle_{\mathcal{H}} = 0 & \text{\rm for $f \in \mathcal{H}_0$,} \\ 
&\langle \mathcal{L} f,f \rangle_{\mathcal{H}} \leq -\frac{1}{n-1} \, \|f\|_{\mathcal{H}}^2 & \text{\rm for $f \in \mathcal{H}_-$.} 
\end{align*}

As in Lemma 2.4 in \cite{Brendle-Choi}, we can show that the function $u(x,\tau)$ satisfies 
\[\frac{\partial}{\partial \tau} u = \mathcal{L} u + E + \langle A(\tau)x,\nu_\Sigma \rangle,\] 
where $E$ is an error term satisfying $|E| \leq O(\rho(\tau)^{-1}) \, (|u| + |\nabla^\Sigma u| + |A(\tau)|)$. We next define $\hat{u}(x,\tau) = u(x,\tau) \, \varphi \big ( \frac{x_{n+1}}{\rho(\tau)} \big )$. The function $\hat{u}(x,\tau)$ satisfies 
\[\frac{\partial}{\partial \tau} \hat{u} = \mathcal{L} \hat{u} + \hat{E} + \langle A(\tau)x,\nu_\Sigma \rangle \, \varphi \Big ( \frac{x_{n+1}}{\rho(\tau)} \Big ),\] 
where $\hat{E}$ is an error term satisfying $\|\hat{E}\|_{\mathcal{H}} \leq O(\rho(\tau)^{-1}) \, (\|\hat{u}\|_{\mathcal{H}} + |A(\tau)|)$ (cf. \cite{Brendle-Choi}, Lemma 2.5). Moreover, the orthogonality relations in Proposition \ref{precise.choice.of.rotation.matrix} imply that the function $\hat{u}(x,\tau)$ is orthogonal (with respect to the inner product on $\mathcal{H}$) to the function $\langle Ax,\nu_\Sigma \rangle$ for every $\tau$ and every matrix $A \in so(n+1)$. \\

\begin{lemma}
\label{pde.for.hat.u.2}
We have $|A(\tau)| \leq O(\rho(\tau)^{-1}) \, \|\hat{u}\|_{\mathcal{H}}$ and 
\[\Big \| \frac{\partial}{\partial \tau} \hat{u} - \mathcal{L} \hat{u} \Big \|_{\mathcal{H}} \leq O(\rho(\tau)^{-1}) \, \|\hat{u}\|_{\mathcal{H}}.\] 
\end{lemma}

\textbf{Proof.} Analogous to \cite{Brendle-Choi}, Lemma 2.6. \\

We now define 
\begin{align*} 
&U_+(\tau) := \|P_+ \hat{u}(\cdot,\tau)\|_{\mathcal{H}}^2, \\ 
&U_0(\tau) := \|P_0 \hat{u}(\cdot,\tau)\|_{\mathcal{H}}^2, \\ 
&U_-(\tau) := \|P_- \hat{u}(\cdot,\tau)\|_{\mathcal{H}}^2, 
\end{align*} 
where $P_+, P_0, P_-$ denote the orthogonal projections to $\mathcal{H}_+,\mathcal{H}_0,\mathcal{H}_-$, respectively. Using Lemma \ref{pde.for.hat.u.2}, we obtain
\begin{align*} 
&\frac{d}{d\tau} U_+(\tau) \geq U_+(\tau) - O(\rho(\tau)^{-1}) \, (U_+(\tau) + U_0(\tau) + U_-(\tau)), \\ 
&\Big | \frac{d}{d\tau} U_0(\tau) \Big | \leq O(\rho(\tau)^{-1}) \, (U_+(\tau) + U_0(\tau) + U_-(\tau)), \\ 
&\frac{d}{d\tau} U_-(\tau) \leq -\frac{2}{n-1} \, U_-(\tau) + O(\rho(\tau)^{-1}) \, (U_+(\tau) + U_0(\tau) + U_-(\tau)). 
\end{align*}
Clearly, $U_+(\tau)+U_0(\tau)+U_-(\tau) = \|\hat{u}\|_{\mathcal{H}}^2 \to 0$ as $\tau \to -\infty$. 

\begin{lemma}
\label{consequence.of.merle.zaag}
We have $U_0(\tau)+U_-(\tau) \leq o(1) U_+(\tau)$.
\end{lemma}

\textbf{Proof.} 
The ODE lemma of Merle and Zaag (cf. Lemma 5.4 in \cite{Angenent-Daskalopoulos-Sesum} or Lemma A.1 in \cite{Merle-Zaag}) implies that either $U_0(\tau)+U_-(\tau) \leq o(1) U_+(\tau)$ or $U_+(\tau)+U_-(\tau) \leq o(1) U_0(\tau)$. 

We now explain how to rule out the second case. If $U_+(\tau)+U_-(\tau) \leq o(1) U_0(\tau)$, then $\frac{\hat{u}(\cdot,\tau)}{\|\hat{u}(\cdot,\tau)\|_{\mathcal{H}}}$ converges with respect to $\|\cdot\|_{\mathcal{H}}$ to the subspace $\mathcal{H}_0 = \text{\rm span}\{x_{n+1}^2-2,x_1x_{n+1},\hdots,x_nx_{n+1}\}$. The orthogonality relations in Proposition \ref{precise.choice.of.rotation.matrix} imply that $\hat{u}(\cdot,\tau)$ is orthogonal to $\langle Ax,\nu_\Sigma \rangle$ for each $A \in so(n+1)$. In other words, $\hat{u}(\cdot,\tau)$ is orthogonal to $x_1x_{n+1},\hdots,x_nx_{n+1}$. Therefore, $\frac{\hat{u}(\cdot,\tau)}{\|\hat{u}(\cdot,\tau)\|_{\mathcal{H}}}$ converges (with respect to $\|\cdot\|_{\mathcal{H}}$) to a non-zero multiple of $x_{n+1}^2-2$.

Let $\tilde{\Omega}_\tau$ denote the region enclosed by $\tilde{M}_\tau$, and let $\mathcal{A}(z,\tau)$ denote the area of the intersection $\tilde{\Omega}_\tau \cap \{x_{n+1}=z\}$. By the Brunn-Minkowski inequality, the function $z \mapsto \mathcal{A}(z,\tau)^{\frac{1}{n}}$ is concave. Since $\tilde{M}_\tau$ is noncompact, it follows that the function $z \mapsto \mathcal{A}(z,\tau)^{\frac{1}{n}}$ is monotone for each $\tau$. 

Note that $\mathcal{A}(z,\tau) = \frac{1}{n} \int_{S^{n-1}} (\sqrt{2(n-1)}+u(\cdot,\tau))^n$ for $|z| \leq \rho(\tau)$. Consequently, the function $z \mapsto \int_{S^{n-1}} [(\sqrt{2(n-1)}+u(\cdot,\tau))^n-\sqrt{2(n-1)}^n] \, d\theta$ is monotone. In particular, we either have 
\begin{align*} 
&\int_{-3}^{-1} \bigg ( \int_{S^{n-1}} [(\sqrt{2(n-1)}+u(\cdot,\tau))^n-\sqrt{2(n-1)}^n] \bigg ) \, dz \\ 
&\leq \int_{-1}^1 \bigg ( \int_{S^{n-1}} [(\sqrt{2(n-1)}+u(\cdot,\tau))^n-\sqrt{2(n-1)}^n] \bigg ) \, dz \\ 
&\leq \int_1^3 \bigg ( \int_{S^{n-1}} [(\sqrt{2(n-1)}+u(\cdot,\tau))^n-\sqrt{2(n-1)}^n] \bigg ) \, dz 
\end{align*} 
or 
\begin{align*} 
&\int_{-3}^{-1} \bigg ( \int_{S^{n-1}} [(\sqrt{2(n-1)}+u(\cdot,\tau))^n-\sqrt{2(n-1)}^n] \bigg ) \, dz \\ 
&\geq \int_{-1}^1 \bigg ( \int_{S^{n-1}} [(\sqrt{2(n-1)}+u(\cdot,\tau))^n-\sqrt{2(n-1)}^n] \bigg ) \, dz \\ 
&\geq \int_1^3 \bigg ( \int_{S^{n-1}} [(\sqrt{2(n-1)}+u(\cdot,\tau))^n-\sqrt{2(n-1)}^n] \bigg ) \, dz. 
\end{align*} 
On the other hand, we know that $\sup |u(\cdot,\tau)| \to 0$, and $\frac{\hat{u}(\cdot,\tau)}{\|\hat{u}(\cdot,\tau)\|_{\mathcal{H}}}$ converges (with respect to $\|\cdot\|_{\mathcal{H}}$) to a non-zero multiple of $x_{n+1}^2-2$. Consequently, the function 
\[\frac{(\sqrt{2(n-1)}+\hat{u}(\cdot,\tau))^n-\sqrt{2(n-1)}^n}{\|\hat{u}(\cdot,\tau)\|_{\mathcal{H}}}\] 
converges (with respect to $\|\cdot\|_{\mathcal{H}}$) to a non-zero multiple of $x_{n+1}^2-2$. This is a contradiction. The proof of Lemma \ref{consequence.of.merle.zaag} is now complete. \\

\begin{lemma} 
\label{almost.sharp.asymptotics.for.u}
For each $\varepsilon > 0$, we have $\|u(\cdot,\tau)\|_{C^4(S^{n-1} \times [-10n,10n])} \leq o(e^{\frac{(1-\varepsilon)\tau}{2}})$ and $|A(\tau)| \leq o(e^{\frac{(1-\varepsilon)\tau}{2}})$.
\end{lemma}

\textbf{Proof.} 
Lemma \ref{consequence.of.merle.zaag} gives $U_0(\tau)+U_-(\tau) \leq o(1) U_+(\tau)$. Substituting this back into the ODE for $U_+(\tau)$ gives  
\[\frac{d}{d\tau} U_+(\tau) \geq U_+(\tau) - o(1) \, U_+(\tau).\] 
Consequently, for every $\varepsilon>0$, we have $U_+(\tau) \leq o(e^{(1-\varepsilon)\tau})$. Using the estimate $U_0(\tau)+U_-(\tau) \leq o(1) U_+(\tau)$, we obtain 
\[\|\hat{u}\|_{\mathcal{H}}^2 = U_+(\tau) + U_0(\tau) + U_-(\tau) \leq o(e^{(1-\varepsilon)\tau}).\] 
This implies $|A(\tau)| \leq o(1) \, \|\hat{u}\|_{\mathcal{H}} \leq o(e^{\frac{(1-\varepsilon)\tau}{2}})$. Moreover, standard interpolation inequalities imply $\|u(\cdot,\tau)\|_{C^4(S^{n-1} \times [-10n,10n])} \leq o(e^{\frac{(1-\varepsilon)\tau}{2}})$. This completes the proof of Lemma \ref{almost.sharp.asymptotics.for.u}. \\

It follows from the estimate $|A(\tau)| \leq o(e^{\frac{(1-\varepsilon)\tau}{2}})$ that the limit $\lim_{\tau \to -\infty} S(\tau)$ exists. Without loss of generality, we may assume that $\lim_{\tau \to -\infty} S(\tau) = \text{\rm id}$. Then $|S(\tau)-\text{\rm id}| \leq o(e^{\frac{(1-\varepsilon)\tau}{2}})$. 

\begin{lemma} 
\label{asymptotics.for.surface}
We have 
\[\sup_{\bar{M}_\tau \cap \{|x_{n+1}| \leq e^{-\frac{\tau}{10}}\}} |x_1^2+\hdots+x_n^2-2(n-1)| \leq e^{\frac{\tau}{10}}\] 
if $-\tau$ is sufficiently large.
\end{lemma}

\textbf{Proof.} 
Lemma \ref{almost.sharp.asymptotics.for.u} implies 
\[\sup_{x \in \bar{M}_\tau \cap B_{5n}(0)} |x_1^2+\hdots+x_n^2-2(n-1)| \leq o(e^{\frac{(1-\varepsilon)\tau}{2}}).\] 
In view of the convexity of $\bar{M}_\tau$, it follows that 
\[\sup_{\bar{M}_\tau \cap \{|x_{n+1}| \leq e^{-\frac{\tau}{10}}\}} (x_1^2+\hdots+x_n^2) \leq 2(n-1) + e^{\frac{\tau}{10}}\] 
if $-\tau$ is sufficiently large. Let 
\[\Sigma_a = \{x \in \mathbb{R}^{n+1}: x_1^2+\hdots+x_n^2=u_a(-x_{n+1})^2, \, -a \leq x_{n+1} \leq 0\}\] 
be the self-similar shrinker constructed in \cite{Angenent-Daskalopoulos-Sesum}. By Lemma 4.4 in \cite{Angenent-Daskalopoulos-Sesum}, $u_a(2) \leq \sqrt{2(n-1)}-a^{-2}$. Since $\bar{M}_\tau$ converges to $\Sigma$ in $C_{\text{\rm loc}}^\infty$, the surface $\bar{M}_\tau \cap \{x_{n+1} \leq -2\}$ encloses the surface $\Sigma_a \cap \{x_{n+1} \leq -2\}$ if $-\tau$ is sufficiently large (depending on $a$). On the other hand, the estimate $\inf_{x \in \bar{M}_\tau \cap B_{5n}(0)} (x_1^2+\hdots+x_n^2) \geq 2(n-1) - o(e^{\frac{(1-\varepsilon)\tau}{2}})$ guarantees that the boundary $\bar{M}_\tau \cap \{x_{n+1}=-2\}$ encloses the boundary $\Sigma_a \cap \{x_{n+1}=-2\}$ provided that $-\tau$ is sufficiently large and $a \leq e^{-\frac{(1-\varepsilon)\tau}{4}}$. By the maximum principle, the surface $\bar{M}_\tau \cap \{x_{n+1} \leq -2\}$ encloses $\Sigma_a \cap \{x_{n+1} \leq -2\}$ whenever $-\tau$ is sufficiently large and $a \leq e^{-\frac{(1-\varepsilon)\tau}{4}}$. By Theorem 8.2 in \cite{Angenent-Daskalopoulos-Sesum}, $u_a(y) \geq \sqrt{2(n-1)(1-a^{-2}y^2)}$. This gives
\[\inf_{\bar{M}_\tau \cap \{-e^{-\frac{\tau}{10}} \leq x_{n+1} \leq -2\}} (x_1^2+\hdots+x_n^2) \geq 2(n-1) - e^{\frac{\tau}{10}}\] 
if $-\tau$ is sufficiently large. An analogous argument gives 
\[\inf_{\bar{M}_\tau \cap \{2 \leq x_{n+1} \leq e^{-\frac{\tau}{10}}\}} (x_1^2+\hdots+x_n^2) \geq 2(n-1) - e^{\frac{\tau}{10}}\] 
if $-\tau$ is sufficiently large. Putting these facts together, we conclude that
\[\inf_{\bar{M}_\tau \cap \{|x_{n+1}| \leq e^{-\frac{\tau}{10}}\}} (x_1^2+\hdots+x_n^2) \geq 2(n-1) - e^{\frac{\tau}{10}}\] 
if $-\tau$ is sufficiently large. This completes the proof of Lemma \ref{asymptotics.for.surface}. \\

\begin{lemma}
\label{size.of.necklike.region} 
Let $\varepsilon_0>0$ be given. If $-\tau$ is sufficiently large (depending on $\varepsilon_0$), then every point in $\bar{M}_\tau \cap \{|x_{n+1}| \leq \frac{1}{2} \, e^{-\frac{\tau}{10}}\}$ lies at the center of an $\varepsilon_0$-neck. Moreover, the radius of that neck is uniformly bounded from above and below.
\end{lemma}

\textbf{Proof.} 
By Lemma \ref{asymptotics.for.surface}, the surface $\bar{M}_\tau \cap \{|x_{n+1}| \leq e^{-\frac{\tau}{10}}\}$ lies inside the cylinder $\{x_1^2+\hdots+x_n^2=2(n-1)+e^{\frac{\tau}{10}}\}$, and outside the cylinder $\{x_1^2+\hdots+x_n^2=2(n-1)-e^{\frac{\tau}{10}}\}$. In other words, the surface $M_t \cap \{|x_{n+1}| \leq (-t)^{\frac{1}{2}+\frac{1}{10}}\}$ lies inside the cylinder $\{x_1^2+\hdots+x_n^2=-2(n-1)t+(-t)^{1-\frac{1}{10}}\}$, and outside the cylinder $\{x_1^2+\hdots+x_n^2=-2(n-1)t-(-t)^{1-\frac{1}{10}}\}$. Since $M_t$ bounds a convex domain, it follows that the surface $M_t \cap \{|x_{n+1}| \leq \frac{4}{5} \, (-t)^{\frac{1}{2}+\frac{1}{10}}\}$ can be written as a graph over the cylinder. Classical estimates of Ecker and Huisken \cite{Ecker-Huisken} then imply that $H \leq C \, (-t)^{-\frac{1}{2}}$ on $M_t \cap \{|x_{n+1}| \leq \frac{3}{5} \, (-t)^{\frac{1}{2}+\frac{1}{10}}\}$. On the other hand, using the noncollapsing condition, we obtain $H \geq \frac{1}{C} \, (-t)^{-\frac{1}{2}}$ on $M_t \cap \{|x_{n+1}| \leq \frac{3}{5} \, (-t)^{\frac{1}{2}+\frac{1}{10}}\}$. Consequently, the mean curvature of the rescaled hypersurface $\bar{M}_\tau \cap \{|x_{n+1}| \leq \frac{3}{5} \, e^{-\frac{\tau}{10}}\}$ is uniformly bounded from above and below. Finally, since $\bar{M}_\tau \cap \{|x_{n+1}| \leq e^{-\frac{\tau}{10}}\}$ is $C^0$-close to a cylinder, it follows that every point on $\bar{M}_\tau \cap \{|x_{n+1}| \leq \frac{1}{2} \, e^{-\frac{\tau}{10}}\}$ lies on an $\varepsilon_0$-neck if $-\tau$ is sufficiently large (depending on $\varepsilon_0$). \\

\begin{proposition} 
\label{sharp.asymptotics.for.rescaled.flow}
We have 
\[\sup_{x \in \bar{M}_\tau \cap B_{5n}(0)} |x_1^2+\hdots+x_n^2-2(n-1)| \leq O(e^{\frac{\tau}{2}}).\] 
\end{proposition} 

\textbf{Proof.} 
We repeat the argument above, this time with $\rho(\tau) = e^{-\frac{\tau}{1000}}$. It follows from Lemma \ref{asymptotics.for.surface} that, in the ball $B_{2e^{-\frac{\tau}{1000}}}(0)$, the surface $\bar{M}_\tau$ can be written as a graph over the cylinder, and the height function has $C^4$-norm less than $e^{\frac{\tau}{100}}$. Arguing as above, we can construct a new function $S(\tau)$ taking values in $SO(n+1)$ with the following properties: 
\begin{itemize}
\item In the ball $B_{2e^{-\frac{\tau}{1000}}}(0)$, the rotated surface $\tilde{M}_\tau = S(\tau) \bar{M}_\tau$ can be written as a graph over the cylinder of some function $u(\cdot,\tau)$; that is,
\[\{x + u(x,\tau) \nu_\Sigma(x): x \in \Sigma \cap B_{2e^{-\frac{\tau}{1000}}}(0)\} \subset \tilde{M}_\tau,\] 
where $\|u(\cdot,\tau)\|_{C^4(\Sigma \cap B_{2e^{-\frac{\tau}{1000}}}(0))} \leq O(e^{\frac{\tau}{200}})$. 
\item The function $u(\cdot,\tau)$ satisfies the orthogonality relations 
\[\int_{\Sigma \cap B_{e^{-\frac{\tau}{1000}}}(0)} e^{-\frac{|x|^2}{4}} \, \langle Ax,\nu_\Sigma \rangle \, u(x,\tau) \, \varphi(e^{\frac{\tau}{1000}} \, x_{n+1}) = 0\] 
for all $A \in so(n+1)$.
\end{itemize}
Moreover, $\lim_{\tau \to -\infty} S(\tau) = \text{\rm id}$. The function $\hat{u}(x,\tau) = u(x,\tau) \, \varphi(e^{\frac{\tau}{1000}} \, x_{n+1})$ satisfies 
\[\Big \| \frac{\partial}{\partial \tau} \hat{u} - \mathcal{L} \hat{u} \Big \|_{\mathcal{H}} \leq O(e^{\frac{\tau}{1000}}) \, \|\hat{u}\|_{\mathcal{H}}.\] 
Hence, if we define 
\begin{align*} 
&U_+(\tau) := \|P_+ \hat{u}(\cdot,\tau)\|_{\mathcal{H}}^2, \\ 
&U_0(\tau) := \|P_0 \hat{u}(\cdot,\tau)\|_{\mathcal{H}}^2, \\ 
&U_-(\tau) := \|P_- \hat{u}(\cdot,\tau)\|_{\mathcal{H}}^2, 
\end{align*} 
then 
\begin{align*} 
&\frac{d}{d\tau} U_+(\tau) \geq U_+(\tau) - O(e^{\frac{\tau}{1000}}) \, (U_+(\tau) + U_0(\tau) + U_-(\tau)), \\ 
&\Big | \frac{d}{d\tau} U_0(\tau) \Big | \leq O(e^{\frac{\tau}{1000}}) \, (U_+(\tau) + U_0(\tau) + U_-(\tau)), \\ 
&\frac{d}{d\tau} U_-(\tau) \leq -U_-(\tau) + O(e^{\frac{\tau}{1000}}) \, (U_+(\tau) + U_0(\tau) + U_-(\tau)). 
\end{align*}
The ODE lemma of Merle and Zaag (cf. Lemma 5.4 in \cite{Angenent-Daskalopoulos-Sesum}) implies that either $U_0(\tau)+U_-(\tau) \leq o(1) U_+(\tau)$ or $U_+(\tau)+U_-(\tau) \leq o(1) U_0(\tau)$. As above, the latter case can be ruled out. Therefore, $U_0(\tau)+U_-(\tau) \leq o(1) U_+(\tau)$. This gives 
\[\frac{d}{d\tau} U_+(\tau) \geq U_+(\tau) - O(e^{\frac{\tau}{1000}}) \, U_+(\tau),\]
hence $U_+(\tau) \leq O(e^\tau)$. Thus, $U_0(\tau)+U_-(\tau) \leq o(1) \, U_+(\tau) \leq O(e^\tau)$. Consequently, $\|\hat{u}\|_{\mathcal{H}} \leq O(e^{\frac{\tau}{2}})$. Arguing as in Lemma \ref{pde.for.hat.u.2}, we obtain $|A(\tau)| \leq O(e^{\frac{\tau}{2}})$. Since $\lim_{\tau \to -\infty} S(\tau) = \text{\rm id}$, we conclude that $|S(\tau)-\text{\rm id}| \leq O(e^{\frac{\tau}{2}})$. Finally, $u$ satisfies an equation of the form $\frac{\partial}{\partial \tau} u = \tilde{\mathcal{L}} u + \langle A(\tau)x,\nu_\Sigma \rangle$, where $\tilde{\mathcal{L}}$ is an elliptic operator of second order whose coefficients depend on $u$, $\nabla u$, $\nabla^2 u$, and $A(\tau)$. As $\tau \to -\infty$, the coefficients of $\tilde{\mathcal{L}}$ converge smoothly to the corresponding coefficients of $\mathcal{L}$. Hence, standard interior estimates for parabolic equations imply that $\|u(\cdot,\tau)\|_{C^4(S^{n-1} \times [-10n,10n])} \leq O(e^{\frac{\tau}{2}})$. Combining this estimate with the estimate $|S(\tau)-\text{\rm id}| \leq O(e^{\frac{\tau}{2}})$, we conclude that 
\[\sup_{x \in \bar{M}_\tau \cap B_{5n}(0)} |x_1^2+\hdots+x_n^2-2(n-1)| \leq O(e^{\frac{\tau}{2}}).\]

\section{Lower bound for $H_{\text{\rm max}}(t)$ as $t \to -\infty$}

\label{barrier}

Let $M_t$, $t \in (-\infty,0]$, be a noncompact ancient solution of mean curvature flow in $\mathbb{R}^{n+1}$ which is strictly convex, uniformly two-convex, and noncollapsed. For each $t$, we define $H_{\text{\rm max}}(t)$ to be the supremum of the mean curvature of $M_t$.

\begin{proposition}
\label{H_max.finite}
For each $t$, $H_{\text{\rm max}}(t) < \infty$.
\end{proposition}

\textbf{Proof.} 
Let us fix a time $t$ and a small number $\varepsilon>0$. It follows from Proposition 3.1 in \cite{Haslhofer-Kleiner2} that every point in $M_t$ which lies outside some large compact set must lie at the center of an $\varepsilon$-neck. Hence, if $H_{\text{\rm max}}(t) = \infty$, then the surface $M_t$ contains a sequence of $\varepsilon$-necks with radii converging to $0$, but this cannot happen in a convex hypersurface. \\

\begin{proposition} 
\label{properties.of.H_max}
The function $H_{\text{\rm max}}(t)$ is continuous and monotone increasing in $t$.
\end{proposition} 

\textbf{Proof.} 
We first show that $H_{\text{\rm max}}(t)$ is continuous in $t$. It follows from work of Haslhofer and Kleiner \cite{Haslhofer-Kleiner1},\cite{Haslhofer-Kleiner2} that $|\frac{\partial}{\partial t} H| \leq CH^3$ for some uniform constant $C$. Let us fix a time $\bar{t}$ and a positive number $\varepsilon>0$. By definition of $H_{\text{\rm max}}(\bar{t})$, there exists a point on $M_{\bar{t}}$ where the mean curvature lies in the interval $[(1-\frac{\varepsilon}{2}) \, H_{\text{\rm max}}(\bar{t}),H_{\text{\rm max}}(\bar{t})]$. If $t$ is sufficiently close to $\bar{t}$, then the estimate $|\frac{\partial}{\partial t} H| \leq CH^3$ implies that there exists a point on $M_t$ where the mean curvature lies in the interval $[(1-\varepsilon) \, H_{\text{\rm max}}(\bar{t}),(1+\varepsilon) \, H_{\text{\rm max}}(\bar{t})]$. In particular, $H_{\text{\rm max}}(t) \geq (1-\varepsilon) \, H_{\text{\rm max}}(\bar{t})$ if $t$ is sufficiently close to $\bar{t}$. Suppose next that there is a point on $M_t$ where the mean curvature is equal to $(1+\varepsilon) \, H_{\text{\rm max}}(t)$. If $t$ is sufficiently close to $\bar{t}$, then the estimate $|\frac{\partial}{\partial t} H| \leq CH^3$ implies that there exists a point on $M_{\bar{t}}$ where the mean curvature lies in the interval $[(1+\frac{\varepsilon}{2}) \, H_{\text{\rm max}}(\bar{t}),(1+2\varepsilon) \, H_{\text{\rm max}}(\bar{t})]$, which contradicts the definition of $H_{\text{\rm max}}(\bar{t})$. Hence, if $t$ is sufficiently close to $\bar{t}$, there is no point on $M_t$ where the mean curvature is equal to $(1+\varepsilon) \, H_{\text{\rm max}}(\bar{t})$. On the other hand, since $M_t$ bounds a convex domain which is not a slab, we know that $M_t$ is connected (cf. \cite{Stoker}, Theorem V). Consequently, $H_{\text{\rm max}}(t) \leq (1+\varepsilon) \, H_{\text{\rm max}}(\bar{t})$ if $t$ is sufficiently close to $\bar{t}$.

Thus, $H_{\text{\rm max}}(t)$ is a continuous function of $t$. In particular, $H_{\text{\rm max}}(t)$ is uniformly bounded from above on every compact time interval. Consequently, $H_{\text{\rm max}}(t)$ is monotone increasing in $t$ by Hamilton's Harnack inequality \cite{Hamilton}. This completes the proof of Proposition \ref{properties.of.H_max}. \\

\begin{proposition}
\label{lower.bound.for.Hmax}
We have $\liminf_{t \to -\infty} H_{\text{\rm max}}(t) > 0$.
\end{proposition} 

\textbf{Proof.} 
Proposition \ref{sharp.asymptotics.for.rescaled.flow} implies that 
\[\sup_{x \in (-t)^{-\frac{1}{2}} \, (M_t \cap B_{5n(-t)^{\frac{1}{2}}}(0))} |x_1^2+\hdots+x_n^2-2(n-1)| \leq O((-t)^{-\frac{1}{2}}).\] 
Since $M_t$ has exactly one end, we can assume without loss of generality that $M_t \cap \{x_{n+1} \geq 0\}$ is noncompact and $M_t \cap \{x_{n+1} \leq 0\}$ is compact. There exists a large constant $K$ with the following property: if $-t$ is sufficiently large, then the cross-section 
\[(-t)^{-\frac{1}{2}} \, (M_t \cap \{x_{n+1}=-2 (-t)^{\frac{1}{2}}\})\] 
lies outside the sphere 
\[\{x_1^2+\hdots+x_n^2=(\sqrt{2(n-1)}-K (-t)^{-\frac{1}{2}})^2, \, x_{n+1}=-2\}.\] 
We now recall the self-similar shrinkers constructed in \cite{Angenent-Daskalopoulos-Sesum}. For $a>0$ large, there exists a self-similar shrinker
\[\Sigma_a = \{x \in \mathbb{R}^{n+1}: x_1^2+\hdots+x_n^2=u_a(-x_{n+1})^2, \, -a \leq x_{n+1} \leq 0\}\] 
satisfying $H = \frac{1}{2} \, \langle x,\nu \rangle$. Consequently, the hypersurfaces 
\begin{align*} 
\Sigma_{a,t} 
&:= (-t)^{\frac{1}{2}} \, \Sigma_a + (0,\hdots,0,Ka^2) \\ 
&= \{x \in \mathbb{R}^{n+1}: x_1^2+\hdots+x_n^2=(-t) \, u_a((-x_{n+1}+Ka^2) (-t)^{-\frac{1}{2}} )^2, \\ 
&\hspace{25mm} Ka^2 - a (-t)^{\frac{1}{2}} \leq x_{n+1} \leq Ka^2\} 
\end{align*} 
evolve by mean curvature flow.

As in \cite{Brendle-Choi}, we can use the hypersurfaces $\Sigma_{a,t} \cap \{x_{n+1} \leq -2 (-t)^{\frac{1}{2}}\}$ as barriers. As $t \to -\infty$, the rescaled surfaces $(-t)^{-\frac{1}{2}} \, M_t$ converge in $C_{\text{\rm loc}}^\infty$ to the cylinder $\{x_1^2+\hdots+x_n^2=2(n-1)\}$. Furthermore, the rescaled surfaces $(-t)^{-\frac{1}{2}} \, (\Sigma_{a,t} \cap \{x_{n+1} \leq -2 (-t)^{\frac{1}{2}}\})$ converge to $\Sigma_a \cap \{x_{n+1} \leq -2\}$ as $t \to -\infty$. It follows from Lemma 4.4 in \cite{Angenent-Daskalopoulos-Sesum} that the set $\Sigma_a \cap \{x_{n+1} \leq -2\}$ is a compact subset of $\{x_1^2+\hdots+x_n^2 < 2(n-1)\}$. Consequently, $\Sigma_{a,t} \cap \{x_{n+1} \leq -2 (-t)^{\frac{1}{2}}\}$ lies inside $M_t \cap \{x_{n+1} \leq -2 (-t)^{\frac{1}{2}}\}$ if $-t$ is sufficiently large (depending on $a$). 

By our choice of $K$, the cross-section 
\[(-t)^{-\frac{1}{2}} \, (M_t \cap \{x_{n+1} = -2 (-t)^{\frac{1}{2}}\})\] 
lies outside the sphere 
\[\{x_1^2+\hdots+x_n^2=(\sqrt{2(n-1)}-K (-t)^{-\frac{1}{2}})^2, \, x_{n+1}=-2\}.\] 
Moreover, the cross-section 
\[(-t)^{-\frac{1}{2}} \, (\Sigma_{a,t} \cap \{x_{n+1} = -2 (-t)^{\frac{1}{2}}\})\] 
is a sphere 
\[\{x_1^2+\hdots+x_n^2=u_a(2 + Ka^2 (-t)^{-\frac{1}{2}})^2, \, x_{n+1}=-2\}.\] 
Using Lemma 4.4 in \cite{Angenent-Daskalopoulos-Sesum}, we obtain $u_a(2) \leq \sqrt{2(n-1)}$ and $u_a(2)-u_a(1) \leq -a^{-2}$ if $a$ is sufficiently large. Since the function $u_a$ is concave, we obtain 
\begin{align*} 
u_a(2+Ka^2 (-t)^{-\frac{1}{2}}) 
&\leq u_a(2) + Ka^2 (-t)^{-\frac{1}{2}} \, (u_a(2)-u_a(1)) \\ 
&\leq \sqrt{2(n-1)} - K (-t)^{-\frac{1}{2}} 
\end{align*} 
for $-t \geq 4K^2a^2$. Consequently, the cross-section $\Sigma_{a,t} \cap \{x_{n+1} = -2 (-t)^{\frac{1}{2}}\}$ lies inside the cross-section $M_t \cap \{x_{n+1} = -2 (-t)^{\frac{1}{2}}\}$ whenever $-t \geq 4K^2a^2$ and $a$ is sufficiently large. By the maximum principle, the hypersurface $\Sigma_{a,t} \cap \{x_{n+1} \leq -2 (-t)^{\frac{1}{2}}\}$ lies inside the hypersurface $M_t \cap \{x_{n+1} \leq -2 (-t)^{\frac{1}{2}}\}$ whenever $-t \geq 4K^2a^2$ and $a$ is sufficiently large. For $-t=4K^2a^2$, the tip of $\Sigma_{a,t}$ has distance $a (-t)^{\frac{1}{2}} -Ka^2 = Ka^2 = -\frac{t}{4K}$ from the origin. Consequently, the intersection $M_t \cap \{x_1=\hdots=x_n=0, \, x_{n+1} \leq \frac{t}{4K}\}$ is non-empty if $-t$ is sufficiently large. In particular, $\limsup_{t \to -\infty} H_{\text{\rm max}}(t) > 0$. Since $H_{\text{\rm max}}(t)$ is monotone increasing in $t$, it follows that $\liminf_{t \to -\infty} H_{\text{\rm max}}(t) > 0$.

\section{The neck improvement theorem}

\label{neck.improvement}

In this section, we show that a neck becomes more symmetric under the evolution.

\begin{definition}
\label{vector.field}
Let $\mathcal{K} = \{K_\alpha: 1 \leq \alpha \leq \frac{n(n-1)}{2}\}$ be a collection of vector fields in $\mathbb{R}^{n+1}$. We say that $\mathcal{K}$ is a normalized set of rotation vector fields if there exists an orthonormal basis $\{J_\alpha: 1 \leq \alpha \leq \frac{n(n-1)}{2}\}$ of $so(n) \subset so(n+1)$, a matrix $S \in O(n+1)$ and a point $q \in \mathbb{R}^{n+1}$ such that 
\[K_\alpha(x) = S J_\alpha S^{-1} (x-q).\] 
\end{definition}

Note that we require that the vector fields $K_\alpha$ all have a common axis of rotation, but we do not require that the axis of rotation passes through the origin.

\begin{lemma}
\label{vector.field.comparison}
For each $n$, we can find a large constant $C$ and a small constant $\varepsilon_0>0$ with the following property. Let $M$ be a hypersurface in $\mathbb{R}^{n+1}$ with boundary. Assume that, after suitable rescaling, $M$ is $\varepsilon_0$-close (in the $C^4$-norm) to a cylinder $S^{n-1} \times [-5,5]$ of radius $1$. Suppose that $\mathcal{K}^{(1)} = \{K_\alpha^{(1)}: 1 \leq \alpha \leq \frac{n(n-1)}{2}\}$ and $\mathcal{K}^{(2)} = \{K_\alpha^{(2)}: 1 \leq \alpha \leq \frac{n(n-1)}{2}\}$ are two normalized sets of rotation vector fields with the following properties: 
\begin{itemize}
\item $\max_\alpha |K_\alpha^{(1)}| \, H \leq 10n$ and $\max_\alpha |K_\alpha^{(2)}| \, H \leq 10n$ at the point $\bar{x}$.
\item $\max_\alpha |\langle K_\alpha^{(1)},\nu \rangle| \, H \leq \varepsilon$ and $\max_\alpha |\langle K_\alpha^{(2)},\nu \rangle| \, H \leq \varepsilon$ on $M$.
\end{itemize}
Then 
\[\inf_{\omega \in O(\frac{n(n-1)}{2})} \sup_{B_{100n H(\bar{x})^{-1}}(\bar{x})} \max_\alpha \bigg | K_\alpha^{(1)} - \sum_{\beta=1}^{\frac{n(n-1)}{2}} \omega_{\alpha\beta} K_\beta^{(2)} \bigg | \, H(\bar{x}) \leq C\varepsilon.\] 
\end{lemma}

\textbf{Proof.} 
The proof is by contradiction. If the assertion is false, then we can find a sequence of hypersurfaces $M^{(l)}$, a sequence of points $\bar{x}_l \in M^{(l)}$, a sequence $\mathcal{K}^{(1,l)} = \{K_\alpha^{(1,l)}: 1 \leq \alpha \leq \frac{n(n-1)}{2}\}$, a sequence $\mathcal{K}^{(2,l)} = \{K_\alpha^{(2,l)}: 1 \leq \alpha \leq \frac{n(n-1)}{2}\}$, and a sequence of positive numbers $\varepsilon_l \to 0$ such that the following statements hold:
\begin{itemize}
\item[(i)] The hypersurfaces $M_l$ converge in $C^4$ to the standard cylinder $S^{n-1} \times [-5,5]$ of radius $1$ around the $x_{n+1}$-axis.
\item[(ii)] $\max_\alpha |K_\alpha^{(1,l)}| \leq 100$ and $\max_\alpha |K_\alpha^{(2,l)}| \leq 100$ at the point $\bar{x}_l$.
\item[(iii)] $\max_\alpha |\langle K_\alpha^{(1,l)},\nu \rangle| \leq \varepsilon_l$ and $\max_\alpha |\langle K_\alpha^{(2,l)},\nu \rangle| \leq \varepsilon_l$ on $M^{(l)}$.
\item[(iv)] We have 
\[\varepsilon_l^{-1} \inf_{\omega \in O(\frac{n(n-1)}{2})} \sup_{B_{1000}(\bar{x}_l)} \max_\alpha \bigg | K_\alpha^{(1,l)} - \sum_{\beta=1}^{\frac{n(n-1)}{2}} \omega_{\alpha\beta} K_\beta^{(2,l)} \bigg | \to \infty\] 
as $l \to \infty$.
\end{itemize} 
Let us write 
\[K_\alpha^{(1,l)}(x) = S^{(1,l)} J_\alpha^{(1,l)} (S^{(1,l)})^{-1} (x-q^{(1,l)})\] 
and 
\[K_\alpha^{(2,l)}(x) = S^{(2,l)} J_\alpha^{(2,l)} (S^{(2,l)})^{-1} (x-q^{(2,l)}).\] 
Here, $\{J_\alpha^{(1,l)}: 1 \leq \alpha \leq \frac{n(n-1)}{2}\}$ and $\{J_\alpha^{(2,l)}: 1 \leq \alpha \leq \frac{n(n-1)}{2}\}$ are orthonormal bases of $so(n) \subset so(n+1)$, $S^{(1,l)},S^{(2,l)} \in O(n+1)$, and $q^{(1,l)},q^{(2,l)} \in \mathbb{R}^{n+1}$. 

The condition (ii) implies that the axis of rotation of $\mathcal{K}^{(1,l)}$ has bounded distance from the origin. It is easy to see that the axis of rotation of $\mathcal{K}^{(1,l)}$ converges to the $x_{n+1}$-axis as $l \to \infty$. Thus, $S^{(1,l)} (0,\hdots,0,1) \to \pm (0,\hdots,0,1)$ as $l \to \infty$. Without loss of generality, we may assume that $S^{(1,l)} \to \text{\rm id}$ as $l \to \infty$. Similarly, we may assume that $S^{(2,l)} \to \text{\rm id}$ as $l \to \infty$. Moreover, we can arrange that $q^{(1,l)}$ and $q^{(2,l)}$ are orthogonal to the $x_{n+1}$-axis. (To achieve this, we add a multiple of $S^{(1,l)} (0,\hdots,0,1)$ to $q^{(1,l)}$, and a multiple of $S^{(2,l)} (0,\hdots,0,1)$ to $q^{(2,l)}$. This does not change the vector fields $K_\alpha^{(1,l)}$ and $K_\alpha^{(2,l)}$.) This gives $q^{(1,l)} \to 0$ and $q^{(2,l)} \to 0$ as $l \to \infty$.

Let us write 
\[(S^{(1,l)})^{-1} S^{(2,l)} = \exp(\sigma^{(l)}) \, U^{(l)},\] 
where $U^{(l)} \in O(n) \subset O(n+1)$ (in other words, $U^{(l)}$ is an isometry of $\mathbb{R}^{n+1}$ which fixes the $x_{n+1}$-axis) and $\sigma^{(l)} \in so(n+1)$ is an anti-symmetric matrix satisfying $\sigma_{ij}^{(l)} = 0$ for $i,j \in \{1,\hdots,n\}$. Finally, we can find a matrix $\omega^{(l)} \in O(\frac{n(n-1)}{2})$ such that 
\[\sum_{\beta=1}^{\frac{n(n-1)}{2}} \omega_{\alpha\beta}^{(l)} \, J_\beta^{(2,l)} = (U^{(l)})^{-1} J_\alpha^{(1,l)} U^{(l)}.\] 
This gives 
\begin{align*} 
\sum_{\beta=1}^{\frac{n(n-1)}{2}} \omega_{\alpha\beta}^{(l)} \, K_\beta^{(2,l)}(x) 
&= S^{(2,l)} (U^{(l)})^{-1} J_\alpha^{(1,l)} U^{(l)} (S^{(2,l)})^{-1} (x-q^{(2,l)}) \\ 
&= S^{(1,l)} \exp(\sigma^{(l)}) J_\alpha^{(1,l)} \exp(-\sigma^{(l)}) (S^{(1,l)})^{-1} (x-q^{(2,l)}). 
\end{align*}
Let $\delta_l := |q^{(1,l)}-q^{(2,l)}|+|\sigma^{(l)}|$. Clearly, 
\[\sup_{B_{1000}(\bar{x}_l)} \max_\alpha \bigg | K_\alpha^{(1,l)} - \sum_{\beta=1}^{\frac{n(n-1)}{2}} \omega_{\alpha\beta}^{(l)} \, K_\beta^{(2,l)} \bigg | \leq C(n) \, \delta_l.\] 
Hence, property (iv) implies $\varepsilon_l^{-1} \, \delta_l \to \infty$ as $l \to \infty$. Let 
\[V_\alpha^{(l)} := \delta_l^{-1} \, \Big ( K_\alpha^{(1,l)} - \sum_{\beta=1}^{\frac{n(n-1)}{2}} \omega_{\alpha\beta}^{(l)} \, K_\beta^{(2,l)} \Big ).\] 
As $l \to \infty$, the vector fields $V_\alpha^{(l)}$ converge to a limiting vector field $V_\alpha$ of the form 
\[V_\alpha(x) = -[\sigma,J_\alpha]x - J_\alpha \zeta,\] 
where $J_\alpha := \lim_{l \to \infty} J_\alpha^{(1,l)}$, $\zeta := \lim_{l \to \infty} \delta_l^{-1} \, (q^{(1,l)}-q^{(2,l)})$, and $\sigma := \lim_{l \to \infty} \delta_l^{-1} \, \sigma^{(l)}$. Note that $\{J_\alpha: 1 \leq \alpha \leq \frac{n(n-1)}{2}\}$ is an orthonormal basis of $so(n) \subset so(n+1)$, $\zeta \in \mathbb{R}^{n+1}$ is orthogonal to the $x_{n+1}$-axis, and $\sigma \in so(n+1)$ is an anti-symmetric matrix satisfying $\sigma_{ij} = 0$ for $i,j \in \{1,\hdots,n\}$. Moreover, by definition of $\delta_l$, we have $|\zeta| + |\sigma| = \lim_{l \to \infty} \delta_l^{-1} \, (|q^{(1,l)}-q^{(2,l)}|+|\sigma^{(l)}|) = 1$.

Using property (iii), we obtain 
\begin{align*} 
\sup_{M^{(l)}} \max_\alpha |\langle V_\alpha^{(l)},\nu \rangle| 
&\leq C(n) \, \delta_l^{-1} \sup_{M^{(l)}} \max_\alpha (|\langle K_\alpha^{(1,l)},\nu \rangle|+|\langle K_\alpha^{(2,l)},\nu \rangle|) \\ 
&\leq C(n) \, \delta_l^{-1} \, \varepsilon_l \\ 
&\to 0. 
\end{align*} 
Consequently, the limiting vector fields $V_\alpha$ are tangential to the cylinder $S^{n-1} \times [-5,5]$. From this, we deduce that $\zeta=0$ and $\sigma=0$. This is a contradiction. This completes the proof of Lemma \ref{vector.field.comparison}. \\

Following \cite{Huisken-Sinestrari3}, pp.~189--190, we denote by $\mathcal{P}(\bar{x},\bar{t},r,\tau)$ the set of all points $(x,t)$ in space-time such that $x \in B_{g(\bar{t})}(\bar{x},r)$ and $t \in [\bar{t}-\tau,\bar{t}]$. With this understood, we define $\hat{\mathcal{P}}(\bar{x},\bar{t},L,\theta) = \mathcal{P}(\bar{x},\bar{t},(n-1)L \, H(\bar{x},\bar{t})^{-1},(n-1)^2\theta \, H(\bar{x},\bar{t})^{-2})$. We say that $(\bar{x},\bar{t})$ lies on an $\varepsilon$-neck if the parabolic neighborhood $\hat{\mathcal{P}}(\bar{x},\bar{t},100,100)$ is, after rescaling, $\varepsilon$-close (in the $C^{10}$-norm), to a family of shrinking cylinders. 

\begin{definition}
Let $M_t$ be a solution of mean curvature flow. We say that a point $(\bar{x},\bar{t})$ is $\varepsilon$-symmetric if there exists a normalized set of rotation vector fields $\mathcal{K} = \{K_\alpha: 1 \leq \alpha \leq \frac{n(n-1)}{2}\}$ such that $\max_\alpha |K_\alpha| \, H \leq 10n$ at the point $(\bar{x},\bar{t})$ and $\max_\alpha |\langle K_\alpha,\nu \rangle| \, H \leq \varepsilon$ in the parabolic neighborhood $\hat{\mathcal{P}}(\bar{x},\bar{t},10,100)$.
\end{definition}

Note that the condition that $\max_\alpha |K_\alpha| \, H \leq 10n$ at the point $(\bar{x},\bar{t})$ ensures that the distance of the point $\bar{x}$ from the axis of rotation of $\mathcal{K} = \{K_\alpha: 1 \leq \alpha \leq \frac{n(n-1)}{2}\}$ is at most $C(n) \, H(\bar{x},\bar{t})^{-1}$.

\begin{theorem}[Neck Improvement Theorem]
\label{neck.improvement.theorem}
Given $n$, we can find a large constant $L$ and a small constant $\varepsilon_1$ with the following property. Suppose that $M_t$ is a solution of mean curvature flow. Moreover, suppose that $(\bar{x},\bar{t})$ is a point in space-time with the property that every point in $\hat{\mathcal{P}}(\bar{x},\bar{t},L,L^2)$ lies at the center of an $\varepsilon_1$-neck and is $\varepsilon$-symmetric, where $\varepsilon \leq \varepsilon_1$. Then $(\bar{x},\bar{t})$ is $\frac{\varepsilon}{2}$-symmetric.
\end{theorem}

\textbf{Proof.} 
Without loss of generality, we assume $\bar{t}=-1$ and $H(\bar{x},-1)=\sqrt{\frac{n-1}{2}}$. We will assume throughout that $L$ is sufficiently large depending on $n$, and $\varepsilon_1$ is sufficiently small depending on $L$. Note that, in the parabolic neighborhood $\hat{\mathcal{P}}(\bar{x},\bar{t},L,L^2)$, the solution can be approximated by a family of shrinking cylinders $S^{n-1}(\sqrt{-2(n-1)t}) \times \mathbb{R}$, up to errors which are bounded by $C(L)\varepsilon_1$ in the $C^{10}$-norm.

\textit{Step 1:} Given any point $(x_0,t_0) \in \hat{\mathcal{P}}(\bar{x},-1,L,L^2)$, we can find a normalized set of rotation vector fields $\mathcal{K}^{(x_0,t_0)} = \{K^{(x_0,t_0)}_\alpha: 1 \leq \alpha \leq \frac{n(n-1)}{2}\}$ such that $\max_\alpha |\langle K_\alpha^{(x_0,t_0)},\nu \rangle| \, H \leq \varepsilon$ on the parabolic neighborhood $\hat{\mathcal{P}}(x_0,t_0,10,100)$. Note that the axis of rotation depends on $(x_0,t_0)$. By a repeated application of Lemma \ref{vector.field.comparison}, we obtain 
\[\inf_{\omega \in O(\frac{n(n-1)}{2})} \sup_{B_{10nL}(0)} \max_\alpha \bigg | K_\alpha^{(\bar{x},-1)} - \sum_{\beta=1}^{\frac{n(n-1)}{2}} \omega_{\alpha\beta} K_\beta^{(x_0,t_0)} \bigg | \leq C(L)\varepsilon\] 
for each point $(x_0,t_0) \in \hat{\mathcal{P}}(\bar{x},-1,L,L^2)$. Without loss of generality, we may assume 
\[\sup_{B_{10nL}(0)} \max_\alpha |K_\alpha^{(\bar{x},-1)} - K_\alpha^{(x_0,t_0)}| \leq C(L)\varepsilon\] 
for each point $(x_0,t_0) \in \hat{\mathcal{P}}(\bar{x},-1,L,L^2)$. This implies 
\[\max_\alpha |\langle K_\alpha^{(\bar{x},-1)},\nu \rangle| \leq C(L)\varepsilon\] 
at  each point $(x_0,t_0) \in \hat{\mathcal{P}}(\bar{x},-1,L,L^2)$. 

For abbreviation, we put $\bar{\mathcal{K}} := \mathcal{K}^{(\bar{x},-1)}$. Without loss of generality, we may assume that the axis of rotation of $\bar{\mathcal{K}} = \{\bar{K}_\alpha: 1 \leq \alpha \leq \frac{n(n-1)}{2}\}$ is the $x_{n+1}$-axis; that is, $\bar{K}_\alpha(x) = J_\alpha x$ for some orthonormal basis $\{J_\alpha: 1 \leq \alpha \leq \frac{n(n-1)}{2}\}$ of $so(n) \subset so(n+1)$. Finally, we may assume that the point $\bar{x}$ lies in the hyperplane $\{x_{n+1}=0\}$.

\textit{Step 2:} It will be convenient to write $M_t$ as a graph over the $x_{n+1}$-axis, so that 
\[\Big \{ (r(\theta,z,t) \, \theta,z): \theta \in S^{n-1}, \, z \in \Big [ -\frac{3L}{4},\frac{3L}{4} \Big ] \Big \} \subset M_t.\] 
Note that the function $r(\theta,z,t) - (-2(n-1)t)^{\frac{1}{2}}$ is bounded by $C(L)\varepsilon_1$ in the $C^{10}$-norm. A straightforward computation gives 
\[\nu = \frac{1}{\sqrt{1+r^{-2} \, |\nabla^{S^{n-1}} r|^2 + (\frac{\partial r}{\partial z})^2}} \, \Big ( \theta - r^{-1} \, \nabla^{S^{n-1}} r,-\frac{\partial r}{\partial z} \Big ),\] 
hence 
\[\langle \bar{K}_\alpha,\nu \rangle = -\frac{\langle J_\alpha \theta,\nabla^{S^{n-1}} r \rangle}{\sqrt{1+r^{-2} \, |\nabla^{S^{n-1}} r|^2 + (\frac{\partial r}{\partial z})^2}},\] 
where $\nabla^{S^{n-1}} r$ represents the gradient of the function $r$ with respect to the angular variables. 

We have shown in Step 1 that $\max_\alpha |\langle \bar{K}_\alpha,\nu \rangle| \leq C(L)\varepsilon$ on the parabolic neighborhood $\hat{\mathcal{P}}(\bar{x},-1,L,L^2)$. Moreover, our assumptions imply that $|\nabla^{S^{n-1}} r| + |\frac{\partial r}{\partial z}| \leq C(L)\varepsilon_1$. This gives $|\nabla^{S^{n-1}} r| \leq C(L)\varepsilon$ and 
\[|\langle \bar{K}_\alpha,\nu \rangle + \langle J_\alpha \theta,\nabla^{S^{n-1}} r \rangle| \leq C(L)\varepsilon_1\varepsilon.\] 
Moreover, the identity $\text{\rm div}_{S^{n-1}}(J_\alpha \theta) = 0$ gives $\text{\rm div}_{S^{n-1}}(r(\theta,z,t) \, J_\alpha \theta) = \langle J_\alpha \theta,\nabla^{S^{n-1}} r \rangle$, hence 
\[|\langle \bar{K}_\alpha,\nu \rangle + \text{\rm div}_{S^{n-1}}(r(\theta,z,t) \, J_\alpha \theta)| \leq C(L)\varepsilon_1\varepsilon.\] 

\textit{Step 3:} Let us fix an index $\alpha \in \{1,\hdots,\frac{n(n-1)}{2}\}$. For each point $(x_0,t_0) \in \hat{\mathcal{P}}(\bar{x},-1,L,L^2)$, the vector field $K^{(x_0,t_0)}_\alpha$ satisfies 
\[|\langle K^{(x_0,t_0)}_\alpha,\nu \rangle| \leq C\varepsilon \, (-t_0)^{\frac{1}{2}}\] 
on the parabolic neighborhood $\hat{\mathcal{P}}(x_0,t_0,10,100)$. On the other hand, we have shown in Step 1 that $\sup_{B_{10nL}(0)} \max_\alpha |\bar{K}_\alpha - K_\alpha^{(x_0,t_0)}| \leq C(L)\varepsilon$. Hence, there exist real numbers $a_1,\hdots,a_n,b_1,\hdots,b_n$ (depending on $(x_0,t_0)$) such that 
\[|a_1|+\hdots+|a_n| \leq C(L) \varepsilon,\] 
\[|b_1|+\hdots+|b_n| \leq C(L) \varepsilon,\] 
and 
\[|\langle \bar{K}_\alpha - K^{(x_0,t_0)}_\alpha,\nu \rangle - (a_1\theta_1+\hdots+a_n\theta_n) - (b_1\theta_1+\hdots+b_n\theta_n)z| \leq C(L)\varepsilon_1\varepsilon\] 
on the parabolic neighborhood $\hat{\mathcal{P}}(x_0,t_0,10,100)$. Consequently, the function $u = \langle \bar{K}_\alpha,\nu \rangle$ satisfies 
\[|u(\theta,z,t) - (a_1\theta_1+\hdots+a_n\theta_n) - (b_1\theta_1+\hdots+b_n\theta_n)z| \leq C\varepsilon \, (-t_0)^{\frac{1}{2}} + C(L)\varepsilon_1\varepsilon\] 
on the parabolic neighborhood $\hat{\mathcal{P}}(x_0,t_0,10,100)$. 

\textit{Step 4:} We again fix an index $\alpha \in \{1,\hdots,\frac{n(n-1)}{2}\}$. The function $u = \langle \bar{K}_\alpha,\nu \rangle$ satisfies the evolution equation 
\[\frac{\partial}{\partial t} u = \Delta_{M_t} u + |A|^2 u.\] 
We have shown in Step 1 that $|u| \leq C(L)\varepsilon$ on the parabolic neighborhood $\hat{P}(\bar{x},-1,L,L^2)$. Using standard interior estimates for parabolic equations, we obtain $|\nabla u| + |\nabla^2 u| \leq C(L) \varepsilon$ on the parabolic neighborhood $\hat{\mathcal{P}}(\bar{x},-1,\frac{L}{2},\frac{L^2}{4})$. This implies 
\[\Big | \frac{\partial}{\partial t} u - \frac{\partial^2}{\partial z^2} u - \frac{1}{(-2(n-1)t)} \, \Delta_{S^{n-1}} u - \frac{1}{(-2t)} \, u \Big | \leq C(L)\varepsilon_1\varepsilon\] 
for $z \in [-\frac{L}{4},\frac{L}{4}]$ and $t \in [-\frac{L^2}{16},-1]$. 

We denote by $\tilde{u}$ the solution of the linear equation 
\[\frac{\partial}{\partial t} \tilde{u} = \frac{\partial^2}{\partial z^2} \tilde{u} + \frac{1}{(-2(n-1)t)} \, \Delta_{S^{n-1}} \tilde{u} + \frac{1}{(-2t)} \, \tilde{u}\] 
in the parabolic cylinder $\{z \in [-\frac{L}{4},\frac{L}{4}], \, t \in [-\frac{L^2}{16},-1]\}$ such that $\tilde{u}=u$ on the parabolic boundary $\{|z|=\frac{L}{4}\} \cup \{t=-\frac{L^2}{16}\}$. The maximum principle gives 
\[|u-\tilde{u}| \leq C(L)\varepsilon_1\varepsilon\] 
in the parabolic cylinder $\{z \in [-\frac{L}{4},\frac{L}{4}], \, t \in [-\frac{L^2}{16},-1]\}$. 

\textit{Step 5:} In order to analyze the PDE for $\tilde{u}$, we perform separation of variables. As above, let $Y_m$ be a basis of eigenfunctions of $\Delta_{S^{n-1}}$, and let $\lambda_m$ denote the corresponding eigenvalues. We assume that the eigenfunctions $Y_m$ are normalized so that $\frac{n}{|S^{n-1}|} \int_{S^{n-1}} Y_m(\theta)^2 \, d\theta = 1$. Clearly, $\sup_{\theta \in S^{n-1}} |Y_m(\theta)| \leq C \, \|Y_m\|_{H^{n-1}(S^{n-1})} \leq C \lambda_m^{n-1}$ for $m \geq 1$. Moreover, $\lambda_m \sim m^{\frac{2}{n-1}}$ as $m \to \infty$. We recall that $\lambda_0=0$, $\lambda_1=\hdots=\lambda_n=n-1$, and $\lambda_{n+1}=2n$. Moreover, $Y_0(\theta)=\frac{1}{\sqrt{n}}$, $Y_1(\theta)=\theta_1, \hdots, Y_n(\theta)=\theta_n$, where $\theta_1,\hdots,\theta_n$ denote the Cartesian coordinates of $\theta \in S^{n-1}$. Let us write 
\[\tilde{u}(\theta,z,t) = \sum_{m=0}^\infty v_m(z,t) \, Y_m(\theta),\] 
where 
\[v_m(z,t) = \frac{n}{|S^{n-1}|} \int_{S^{n-1}} \tilde{u}(\theta,z,t) \, Y_m(\theta) \, d\theta.\] 
Then 
\[\frac{\partial}{\partial t} v_m(z,t) = \frac{\partial^2}{\partial z^2} v_m(z,t) + \frac{n-1-\lambda_m}{(-2(n-1)t)} \, v_m(z,t).\] 
Hence, the rescaled function $\hat{v}_m(z,t) = (-t)^{\frac{1}{2}-\frac{\lambda_m}{2(n-1)}} \, v_m(z,t)$ satisfies 
\[\frac{\partial}{\partial t} \hat{v}_m(z,t) = \frac{\partial^2}{\partial z^2} \hat{v}_m(z,t).\] 

We first consider the case when $m \geq n+1$, so that $\lambda_m \geq 2n$. Using the results in Step 3 and Step 4, we obtain 
\[|v_m(z,t)| \leq (C\varepsilon + C(L)\varepsilon_1\varepsilon)  \, (-t)^{\frac{1}{2}},\] 
and consequently 
\[|\hat{v}_m(z,t)| \leq (C\varepsilon + C(L)\varepsilon_1\varepsilon)  \, (-t)^{1-\frac{\lambda_m}{2(n-1)}}\] 
in the parabolic cylinder $\{z \in [-\frac{L}{4},\frac{L}{4}], \, t \in [-\frac{L^2}{16},-1]\}$. Using the solution formula for the one-dimensional heat equation with Dirichlet boundary condition on the rectangle $[-\frac{L}{4},\frac{L}{4}] \times [-\frac{L^2}{16},-1]$, we obtain 
\begin{align*}
|\hat{v}_m(z,t)| 
&\leq C \, L^{-1} \int_{-\frac{L}{4}}^{\frac{L}{4}} |\hat{v}_m(z,-\frac{L^2}{16})| \, dz \\ 
&+ C \, L \int_{-\frac{L^2}{16}}^t e^{-\frac{L^2}{100(t-s)}} \, (t-s)^{-\frac{3}{2}} \, \Big ( |\hat{v}_m(\frac{L}{4},s)| + |\hat{v}_m(-\frac{L}{4},s)| \Big ) \, ds, 
\end{align*} 
hence 
\begin{align*} 
|\hat{v}_m(z,t)|
&\leq (C \varepsilon+C(L) \varepsilon_1 \varepsilon) \, \Big ( \frac{L}{4} \Big )^{2-\frac{\lambda_m}{n-1}}  \\ 
&+ (C\varepsilon+C(L) \varepsilon_1 \varepsilon) L \int_{-\frac{L^2}{16}}^t e^{-\frac{L^2}{100(t-s)}} \, (t-s)^{-\frac{3}{2}} \, (-s)^{1-\frac{\lambda_m}{2(n-1)}} \, ds \\
&\leq (C \varepsilon+C(L) \varepsilon_1 \varepsilon) \, \Big ( \frac{L}{4} \Big )^{2-\frac{\lambda_m}{n-1}} \\ 
&+ (C\varepsilon+C(L) \varepsilon_1 \varepsilon) L^{-\frac{1}{n-1}} \int_{-\frac{L^2}{16}}^t e^{-\frac{L^2}{200(t-s)}} \, (-s)^{\frac{1-\lambda_m}{2(n-1)}} \, ds \\ 
&\leq (C \varepsilon+C(L) \varepsilon_1 \varepsilon) \, \Big ( \frac{L}{4} \Big )^{2-\frac{\lambda_m}{n-1}} \\ 
&+ (C\varepsilon+C(L) \varepsilon_1 \varepsilon) L^{-\frac{1}{n-1}} \int_{-\frac{L^2}{16}}^{(1+\frac{1}{\sqrt{\lambda_m}})t} e^{-\frac{L^2}{200(t-s)}} \, (-s)^{\frac{1-\lambda_m}{2(n-1)}} \, ds \\ 
&+ (C\varepsilon+C(L) \varepsilon_1 \varepsilon) L^{-\frac{1}{n-1}} \int_{(1+\frac{1}{\sqrt{\lambda_m}})t}^t e^{-\frac{L^2}{200(t-s)}} \, (-s)^{\frac{1-\lambda_m}{2(n-1)}} \, ds \\ 
&\leq (C \varepsilon+C(L) \varepsilon_1 \varepsilon) \, \Big ( \frac{L}{4} \Big )^{2-\frac{\lambda_m}{n-1}} \\ 
&+ (C\varepsilon+C(L) \varepsilon_1 \varepsilon) L^{-\frac{1}{n-1}}  \, \Big ( 1+\frac{1}{\sqrt{\lambda_m}} \Big )^{\frac{2n-1-\lambda_m}{2(n-1)}} \, (-t)^{\frac{2n-1-\lambda_m}{2(n-1)}} \\ 
&+ (C\varepsilon+C(L) \varepsilon_1 \varepsilon) L^{-\frac{1}{n-1}} \, e^{-\frac{L^2\sqrt{\lambda_m}}{200(-t)}} \, (-t)^{\frac{2n-1-\lambda_m}{2(n-1)}} 
\end{align*}
for $z \in [-20n,20n]$ and $t \in [-400n^2,-1]$. This implies 
\begin{align*} 
|v_m(z,t)|
&\leq (C \varepsilon+C(L) \varepsilon_1 \varepsilon) \, \Big ( \frac{L^2}{16(-t)} \Big )^{1-\frac{\lambda_m}{2(n-1)}} \\ 
&+ (C\varepsilon+C(L) \varepsilon_1 \varepsilon) L^{-\frac{1}{n-1}}  \, \Big ( 1+\frac{1}{\sqrt{\lambda_m}} \Big )^{\frac{2n-1-\lambda_m}{2(n-1)}} \\ 
&+ (C\varepsilon+C(L) \varepsilon_1 \varepsilon) L^{-\frac{1}{n-1}} \, e^{-\frac{L^2\sqrt{\lambda_m}}{200(-t)}}  
\end{align*}
for $z \in [-20n,20n]$ and $t \in [-400n^2,-1]$. We now sum over all $m \geq n+1$. Using the estimate $\sup_{\theta \in S^{n-1}} |Y_m(\theta)| \leq C \, \lambda_m^{n-1}$, we obtain 
\[\bigg | \sum_{m=n+1}^\infty v_m(z,t) \, Y_m(\theta) \bigg | \leq C \sum_{m=n+1}^\infty \lambda_m^{n-1} \, |v_m(z,t)| \leq CL^{-\frac{1}{n-1}} \varepsilon + C(L)\varepsilon_1\varepsilon\]
for $z \in [-20n,20n]$ and $t \in [-400n^2,-1]$.

We next consider the case when $1 \leq m \leq n$, so that $\lambda_m=n-1$ and $Y_m(\theta) = \theta_m$. In this case, the function $v_m(z,t)$ satisfies 
\[\frac{\partial}{\partial t} v_m(z,t) = \frac{\partial^2}{\partial z^2} v_m(z,t).\] 
The results in Step 3 and Step 4 imply that, given any point $(z_0,t_0) \in [-\frac{L}{4},\frac{L}{4}] \times [-\frac{L^2}{16},-1]$, we can find real numbers $a_1,\hdots,a_n,b_1,\hdots,b_n$ (depending on $(z_0,t_0)$) such that 
\[|a_1|+\hdots+|a_n| \leq C(L) \varepsilon,\] 
\[|b_1|+\hdots+|b_n| \leq C(L) \varepsilon,\] 
and 
\[|v_m(z,t) - (a_m+b_mz)| \leq C\varepsilon (-t_0)^{\frac{1}{2}} + C(L)\varepsilon_1\varepsilon\] 
for $z \in [z_0-(-t_0)^{\frac{1}{2}},z_0+(-t_0)^{\frac{1}{2}}]$ and $t \in [2t_0,t_0]$. Using interior estimates for the linear heat equation, we obtain 
\[\Big | \frac{\partial^2 v_m}{\partial z^2}(z,t) \Big | \leq (C\varepsilon + C(L)\varepsilon_1\varepsilon) \, (-t)^{-\frac{1}{2}}\] 
in the parabolic cylinder $\{z \in [-\frac{L}{4},\frac{L}{4}], \, t \in [-\frac{L^2}{16},-1]\}$. Using the solution formula for the one-dimensional heat equation with Dirichlet boundary condition on the rectangle $[-\frac{L}{4},\frac{L}{4}] \times [-\frac{L^2}{16},-1]$, we obtain 
\begin{align*} 
\Big | \frac{\partial^2 v_m}{\partial z^2}(z,t) \Big | 
&\leq (C\varepsilon + C(L)\varepsilon_1\varepsilon) \, \Big ( \frac{L}{4} \Big )^{-1} \\ 
&+ (C\varepsilon+C(L) \varepsilon_1 \varepsilon) L \int_{-\frac{L^2}{16}}^t e^{-\frac{L^2}{100(t-s)}} \, (t-s)^{-\frac{3}{2}} \, (-s)^{-\frac{1}{2}} \, ds \\ 
&\leq (C\varepsilon + C(L)\varepsilon_1\varepsilon) \, \Big ( \frac{L}{4} \Big )^{-1} \\ 
&+ (C\varepsilon+C(L) \varepsilon_1 \varepsilon) L^{-2} \int_{-\frac{L^2}{16}}^t (-s)^{-\frac{1}{2}} \, ds \\ 
&\leq (C\varepsilon + C(L)\varepsilon_1\varepsilon) \, L^{-1} 
\end{align*} 
for $z \in [-20n,20n]$ and $t \in [-400n^2,-1]$. Consequently, there exist real numbers $A_1,\hdots,A_n,B_1,\hdots,B_n$ such that 
\[|v_m(z,t) - (A_m+B_m z)| \leq CL^{-1} \varepsilon + C(L)\varepsilon_1\varepsilon\] 
for $z \in [-20n,20n]$ and $t \in [-400n^2,-1]$. 

Finally, we consider the case $m=0$. By the results in Step 2, the function $u = \langle \bar{K}_\alpha,\nu \rangle$ satisfies 
\[|u(\theta,z,t) + \text{\rm div}_{S^{n-1}}(r(\theta,z,t) \, J_\alpha \theta)| \leq C(L)\varepsilon_1\varepsilon.\] 
Integrating over $\theta \in S^{n-1}$ gives 
\[\bigg | \int_{S^{n-1}} u(\theta,z,t) \, d\theta \bigg | \leq C(L)\varepsilon_1\varepsilon,\] 
hence 
\[|v_0(z,t)| \leq C(L)\varepsilon_1\varepsilon\] 
for $z \in [-20n,20n]$ and $t \in [-400n^2,-1]$. 

Putting everything together, we conclude that 
\[|\tilde{u}(\theta,z,t) - (A_1 \theta_1+\hdots+A_n \theta_n) - (B_1 \theta_1+\hdots+B_n \theta_n)z| \leq CL^{-\frac{1}{n-1}} \varepsilon + C(L)\varepsilon_1\varepsilon\] 
for $z \in [-20n,20n]$ and $t \in [-400n^2,-1]$.

\textit{Step 6:} By combining the results in Step 4 and Step 5, we can draw the following conclusion. For each $\alpha \in \{1,\hdots,\frac{n(n-1)}{2}\}$ we can find real numbers $A_{\alpha,1},\hdots,A_{\alpha,n},B_{\alpha,1},\hdots,B_{\alpha,n}$ such that 
\[|A_{\alpha,1}|+\hdots+|A_{\alpha,n}| \leq C(L)\varepsilon,\] 
\[|B_{\alpha,1}|+\hdots+|B_{\alpha,n}| \leq C(L)\varepsilon,\] 
and 
\begin{align*} 
&|u_\alpha(\theta,z,t) - (A_{\alpha,1}\theta_1+\hdots+A_{\alpha,n}\theta_n) - (B_{\alpha,1}\theta_1+\hdots+B_{\alpha,n}\theta_n)z| \\ 
&\leq CL^{-\frac{1}{n-1}} \varepsilon + C(L)\varepsilon_1\varepsilon 
\end{align*}
for $z \in [-20n,20n]$ and $t \in [-400n^2,-1]$, where $u_\alpha(\theta,z,t) := \langle \bar{K}_\alpha,\nu \rangle$. 

\textit{Step 7:} For each $i \in \{1,\hdots,n\}$, we define 
\[E_i = \frac{n}{|S^{n-1}|} \int_{S^{n-1}} r(\theta,0,-1) \, \theta_i \, d\theta\] 
and 
\[F_i = \frac{n}{2 \, |S^{n-1}|} \int_{S^{n-1}} [r(\theta,1,-1)-r(\theta,-1,-1)] \, \theta_i \, d\theta.\] 
We have shown in Step 2 that $|\nabla^{S^{n-1}} r| \leq C(L) \varepsilon$. This implies $|E_i|,|F_i| \leq C(L)\varepsilon$ for $i \in \{1,\hdots,n\}$. By the results in Step 2, the function $u_\alpha = \langle \bar{K}_\alpha,\nu \rangle$ satisfies 
\[|u_\alpha(\theta,z,t) + \text{\rm div}_{S^{n-1}}(r(\theta,z,t) \, J_\alpha \theta)| \leq C(L)\varepsilon_1\varepsilon.\] 
A direct calculation gives 
\[\text{\rm div}_{S^{n-1}}(r(\theta,z,t) \, J_\alpha \theta) \, \theta_i = \text{\rm div}_{S^{n-1}}(r(\theta,z,t) \, \theta_i \, J_\alpha \theta) - r(\theta,z,t) \sum_{j=1}^n J_{\alpha,ij} \, \theta_j,\] 
where $J_{\alpha,ij}$ denote the components of the anti-symmetric matrix $J_\alpha$. Putting these facts together, we obtain 
\[\Big | u_\alpha(\theta,z,t) \, \theta_i + \text{\rm div}_{S^{n-1}}(r(\theta,z,t) \, \theta_i \, J_\alpha \theta) - r(\theta,z,t) \sum_{j=1}^n J_{\alpha,ij} \, \theta_j \Big | \leq C(L)\varepsilon_1\varepsilon\] 
for all $i \in \{1,\hdots,n\}$. Integrating over $\theta \in S^{n-1}$ gives
\[\max_{\alpha,i} \bigg | \frac{n}{|S^{n-1}|} \int_{S^{n-1}} u_\alpha(\theta,0,-1) \, \theta_i \, d\theta - \sum_{j=1}^n J_{\alpha,ij} E_j \bigg | \leq C(L) \varepsilon_1 \varepsilon\] 
and 
\[\max_{\alpha,i} \bigg | \frac{n}{2 \, |S^{n-1}|} \int_{S^{n-1}} [u_\alpha(\theta,1,-1)-u_\alpha(\theta,-1,-1)] \, \theta_i \, d\theta - \sum_{j=1}^n J_{\alpha,ij} F_j \bigg | \leq C(L) \varepsilon_1 \varepsilon.\] 
On the other hand, using the estimate for $u_\alpha(\theta,z,t) - (A_{\alpha,1}\theta_1+\hdots+A_{\alpha,n}\theta_n) - (B_{\alpha,1}\theta_1+\hdots+B_{\alpha,n}\theta_n)z$ in Step 6, we obtain 
\[\max_{\alpha,i} \bigg | \frac{n}{|S^{n-1}|} \int_{S^{n-1}} u_\alpha(\theta,0,-1) \, \theta_i \, d\theta - A_{\alpha,i} \bigg | \leq CL^{-\frac{1}{n-1}} \varepsilon + C(L)\varepsilon_1\varepsilon\] 
and 
\begin{align*} 
&\max_{\alpha,i} \bigg | \frac{n}{2 \, |S^{n-1}|} \int_{S^{n-1}} [u_\alpha(\theta,1,-1)-u_\alpha(\theta,-1,-1)] \, \theta_i \, d\theta - B_{\alpha,i} \bigg | \\ 
&\leq CL^{-\frac{1}{n-1}} \varepsilon + C(L)\varepsilon_1\varepsilon. 
\end{align*}
Putting these facts together, we obtain
\[\max_{\alpha,i} \Big | A_{\alpha,i} - \sum_{j=1}^n J_{\alpha,ij} E_j \Big | \leq CL^{-\frac{1}{n-1}} \varepsilon + C(L)\varepsilon_1\varepsilon\] 
and 
\[\max_{\alpha,i} \Big | B_{\alpha,i} - \sum_{j=1}^n J_{\alpha,ij} F_j \Big | \leq CL^{-\frac{1}{n-1}} \varepsilon + C(L)\varepsilon_1\varepsilon.\] 
Substituting this back into the estimate for $u_\alpha(\theta,z,t) - (A_{\alpha,1}\theta_1+\hdots+A_{\alpha,n}\theta_n) - (B_{\alpha,1}\theta_1+\hdots+B_{\alpha,n}\theta_n)z$ in Step 6, we finally conclude 
\[\max_\alpha \Big | \langle \bar{K}_\alpha,\nu \rangle - \sum_{i,j=1}^n J_{\alpha,ij}  (E_j \theta_i + F_j \theta_iz) \Big | \leq CL^{-\frac{1}{n-1}} \varepsilon + C(L)\varepsilon_1\varepsilon\]
for $z \in [-20n,20n]$ and $t \in [-400n^2,-1]$.

\textit{Step 8:} Let us define a vector $q \in \mathbb{R}^{n+1}$ by $q_i = E_i$ for $i \in \{1,\hdots,n\}$ and $q_{n+1} = 0$. Moreover, let us define an anti-symmetric matrix $\sigma \in so(n+1)$ by $\sigma_{ij} = 0$ for $i,j \in \{1,\hdots,n\}$ and $\sigma_{i,n+1} = F_i$ for $i \in \{1,\hdots,n\}$. Clearly, $|q| \leq C(L)\varepsilon$ and $|\sigma| \leq C(L)\varepsilon$. We define a normalized set of rotation vector fields $\tilde{\mathcal{K}} = \{\tilde{K}_\alpha: 1 \leq \alpha \leq \frac{n(n-1)}{2}\}$ by 
\[\tilde{K}_\alpha(x) := SJ_\alpha S^{-1} (x-q),\] 
where $S := \exp(\sigma) \in SO(n+1)$. Then 
\[\max_\alpha \Big | \langle \tilde{K}_\alpha - \bar{K}_\alpha,\nu \rangle + \sum_{i,j=1}^n J_{\alpha,ij}  (E_j \theta_i + F_j \theta_iz) \Big | \leq C(L)\varepsilon_1\varepsilon\] 
for $z \in [-20n,20n]$ and $t \in [-400n^2,-1]$. Combining this estimate with the estimate in Step 7, we conclude that 
\[\max_\alpha |\langle \tilde{K}_\alpha,\nu \rangle| \leq CL^{-\frac{1}{n-1}} \varepsilon + C(L)\varepsilon_1\varepsilon\] 
for $z \in [-20n,20n]$ and $t \in [-400n^2,-1]$. Consequently, the point $(\bar{x},-1)$ is $(CL^{-\frac{1}{n-1}} \varepsilon + C(L) \varepsilon_1 \varepsilon)$-symmetric. In particular, if we choose $L$ sufficiently large and $\varepsilon_1$ sufficiently small (depending on $L$), then $(\bar{x},-1)$ is $\frac{\varepsilon}{2}$-symmetric. This completes the proof of the Neck Improvement Theorem.

\section{Proof of rotational symmetry}

\label{rotational.symmetry}

Let $M_t$, $t \in (-\infty,0]$, be a noncompact ancient solution of mean curvature flow in $\mathbb{R}^{n+1}$ which is strictly convex, uniformly two-convex, and noncollapsed. For each $t$, the hypersurface $M_t$ bounds a convex domain which we denote by $\Omega_t$. As in \cite{Brendle-Choi}, if $-t$ is sufficiently large, there exists a unique point $p_t \in M_t$ where the mean curvature attains its maximum. Moreover, this is a non-degenerate maximum in the sense that the Hessian of the mean curvature at $p_t$ is negative definite, and the hypersurface $M_t$ looks like the bowl soliton near $p_t$. 

Let $\varepsilon_1$ and $L$ be the constants in the Neck Improvement Theorem. Recall that $H_{\text{\rm max}}(t)$ is uniformly bounded from below. By Proposition 3.1 in \cite{Haslhofer-Kleiner2}, we can find a large constant $\Lambda$ such that the following holds. If $x$ is a point on $M_t$ such that $|x-p_t| \geq \Lambda$, then $x$ lies at the center of an $\varepsilon_1$-neck and furthermore $H(x,t) \, |x-p_t| \geq 10^6 \, n^2 L$. 

\begin{lemma}
\label{derivative.of.distance}
There exists a time $T<0$ with the property that $\frac{d}{dt} |x-p_t| < 0$ whenever $t \leq T$, $x \in M_t \cup \Omega_t$, and $|x-p_t| \geq \Lambda$.
\end{lemma}

\textbf{Proof.} 
If $-t$ is sufficiently large, then $M_t$ looks like the bowl soliton near the point $p_t$. Hence, if $-t$ is sufficiently large, then the vector $\frac{d}{dt} p_t$ is almost parallel to $-\nu(p_t,t)$. Moreover, using the convexity of $\Omega_t$, we obtain $-\langle x-p_t,\nu(p_t,t) \rangle \geq c \, |x-p_t|$ whenever $t \leq T$, $x \in M_t \cup \Omega_t$, and $|x-p_t| \geq \Lambda$. Here, $c$ is a small positive constant. Putting these facts together, we conclude that $\langle x-p_t,\frac{d}{dt} p_t \rangle > 0$ whenever $t \leq T$, $x \in M_t \cup \Omega_t$, and $|x-p_t| \geq \Lambda$. This gives $\frac{1}{2} \frac{d}{dt} |x-p_t|^2 = -\langle x-p_t,\frac{d}{dt} p_t \rangle < 0$ whenever $t \leq T$, $x \in M_t \cup \Omega_t$, and $|x-p_t| \geq \Lambda$. \\

\begin{lemma}
\label{distance.to.tip.decreases}
Let $T$ be defined as in Lemma \ref{derivative.of.distance}. Suppose that $\bar{t} \leq T$, and $\bar{x}$ is a point on $M_{\bar{t}}$ satisfying $|\bar{x}-p_{\bar{t}}| \geq \Lambda$. Then $|\bar{x}-p_t| \geq |\bar{x}-p_{\bar{t}}|$ for all $t \leq \bar{t}$.
\end{lemma} 

\textbf{Proof.} 
We argue by contradiction. Suppose that the assertion is false, and let $\tilde{t} := \sup \{t \leq \bar{t}: |\bar{x}-p_t| < |\bar{x}-p_{\bar{t}}|\}$. By Lemma \ref{derivative.of.distance}, we have $\frac{d}{dt} |\bar{x}-p_t| \big |_{t=\bar{t}} < 0$. Consequently, $\tilde{t} < \bar{t}$. Moreover, for each $t \in [\tilde{t},\bar{t}]$, we have $\bar{x} \in M_t \cup \Omega_t$ and $|\bar{x}-p_t| \geq |\bar{x}-p_{\bar{t}}| \geq \Lambda$. Using Lemma \ref{derivative.of.distance}, we obtain $\frac{d}{dt} |\bar{x}-p_t| < 0$ for all $t \in [\tilde{t},\bar{t}]$. Consequently, $|\bar{x}-p_{\tilde{t}}| > |\bar{x}-p_{\bar{t}}|$. This contradicts the definition of $\tilde{t}$. This completes the proof of Lemma \ref{distance.to.tip.decreases}. \\

\begin{proposition} 
\label{iteration}
Let $T$ be defined as in Lemma \ref{derivative.of.distance}. If $t \leq T$, $x \in M_t$ and $|x-p_t| \geq 2^{\frac{j}{400}} \, \Lambda$, then $(x,t)$ is $2^{-j} \varepsilon_1$-symmetric. 
\end{proposition}

\textbf{Proof.} 
The proof is by induction on $j$. It follows from our choice of $\Lambda$ that the assertion is true for $j=0$. Suppose now that $j \geq 1$ and the assertion holds for $j-1$. We claim that the assertion holds for $j$. Suppose this is false. Then there exists a time $\bar{t} \leq T$ and a point $\bar{x} \in M_{\bar{t}}$ such that $|\bar{x}-p_{\bar{t}}| \geq 2^{\frac{j}{400}} \, \Lambda$ and $(\bar{x},\bar{t})$ is not $2^{-j} \varepsilon_1$-symmetric. By the Neck Improvement Theorem, there exists a point $(x,t) \in \hat{\mathcal{P}}(\bar{x},\bar{t},L,L^2)$ such that either $(x,t)$ is not $2^{-j+1} \varepsilon_1$-symmetric or $(x,t)$ does not lie at the center of an $\varepsilon_1$-neck. In view of the induction hypothesis, we conclude that $|x-p_t| \leq 2^{\frac{j-1}{400}} \, \Lambda$. Since $t \leq \bar{t} \leq T$, Lemma \ref{distance.to.tip.decreases} gives $|\bar{x}-p_{\bar{t}}| \leq |\bar{x}-p_t|$. Putting these facts together, we obtain 
\begin{align*} 
|\bar{x}-p_{\bar{t}}| 
&\leq |\bar{x}-p_t| \\ 
&\leq |x-p_t|+|x-\bar{x}| \\ 
&\leq 2^{\frac{j-1}{400}} \, \Lambda + 10 \, n^2L \, H(\bar{x},\bar{t})^{-1} \\ 
&\leq 2^{-\frac{1}{400}} \, |\bar{x}-p_{\bar{t}}| + \frac{1}{1000} \, |\bar{x}-p_{\bar{t}}| \\ 
&< |\bar{x}-p_{\bar{t}}|. 
\end{align*}
This is a contradiction. \\

\begin{theorem} 
\label{symmetry}
Let $T$ be defined as in Lemma \ref{derivative.of.distance}. Then the hypersurface $M_t$ is rotationally symmetric for each $t \leq T$.
\end{theorem}

\textbf{Proof.} 
The argument is similar to \cite{Brendle-Choi}. We fix a time $\bar{t} \leq T$. For each $j$, let $\Omega^{(j)}$ be the set of all points $(x,t)$ in space-time satisfying $t \leq \bar{t}$ and $|x-p_t| \leq 2^{\frac{j}{400}} \, \Lambda$. If $j$ is sufficiently large, then $H(x,t) \geq n \cdot 2^{-\frac{j}{400}}$ for each point $(x,t) \in \Omega^{(j)}$. Proposition \ref{iteration} guarantees that every point $(x,t) \in \partial \Omega^{(j)}$ is $2^{-j} \varepsilon_1$-symmetric. Consequently, given any point $(x,t) \in \partial \Omega^{(j)}$, we can find a normalized set of rotation vector fields $\mathcal{K}^{(x,t)} = \{K_\alpha^{(x,t)}: 1 \leq \alpha \leq \frac{n(n-1)}{2}\}$ such that $\max_\alpha |\langle K_\alpha^{(x,t)},\nu \rangle| \, H \leq 2^{-j} \varepsilon_1$ on $\hat{\mathcal{P}}(x,t,10,100)$. Using Lemma \ref{vector.field.comparison}, we are able to control how the axis of rotation of $\mathcal{K}^{(x,t)}$ varies as we vary the point $(x,t)$. More precisely, if $(x_1,t_1)$ and $(x_2,t_2)$ are points in $\partial \Omega^{(j)}$ satisfying $(x_2,t_2) \in \hat{\mathcal{P}}(x_1,t_1,\frac{1}{10n},\frac{1}{100n^2})$, then 
\begin{align*} 
&\inf_{\omega \in O(\frac{n(n-1)}{2})} \sup_{B_{10n H(x_2,t_2)^{-1}}(x_2)} \max_\alpha \bigg | K_\alpha^{(x_1,t_1)} - \sum_{\beta=1}^{\frac{n(n-1)}{2}} \omega_{\alpha\beta} K_\beta^{(x_2,t_2)} \bigg | \\ 
&\leq C \, 2^{-j} \, H(x_2,t_2)^{-1}. 
\end{align*}
Now, if $(x_2,t_2)$ is a point in $\partial \Omega^{(j)}$, then the mean curvature $H(x_2,t_2)$ satisfies the estimate $\frac{1}{C} \, 2^{-\frac{j}{400}} \leq H(x_2,t_2) \leq C$. Hence, if $(x_1,t_1)$ and $(x_2,t_2)$ are points in $\partial \Omega^{(j)}$ satisfying $(x_2,t_2) \in \hat{\mathcal{P}}(x_1,t_1,\frac{1}{10n},\frac{1}{100n^2})$, then 
\begin{align*} 
&\inf_{\omega \in O(\frac{n(n-1)}{2})} \sup_{B_{1/C}(x_2)} \max_\alpha \bigg | K_\alpha^{(x_1,t_1)} - \sum_{\beta=1}^{\frac{n(n-1)}{2}} \omega_{\alpha\beta} K_\beta^{(x_2,t_2)} \bigg | \\ 
&\leq C \, 2^{-j+\frac{j}{400}},
\end{align*} 
and consequently 
\begin{align*} 
&\inf_{\omega \in O(\frac{n(n-1)}{2})} \sup_{B_{2^{j/20}}(x_2)} \max_\alpha \bigg | K_\alpha^{(x_1,t_1)} - \sum_{\beta=1}^{\frac{n(n-1)}{2}} \omega_{\alpha\beta} K_\beta^{(x_2,t_2)} \bigg | \\ 
&\leq C \, 2^{-j+\frac{j}{10}}.
\end{align*} 
By a repeated application of this estimate, we can show that, for each $j$, there exists a normalized set of rotation vector fields $\mathcal{K}^{(j)} = \{K_\alpha^{(j)}: 1 \leq \alpha \leq \frac{n(n-1)}{2}\}$ with the following property: if $(x,t)$ is a point in $\partial \Omega^{(j)}$ satisfying $\bar{t}-2^{\frac{j}{100}} \leq t \leq \bar{t}$, then 
\[\inf_{\omega \in O(\frac{n(n-1)}{2})} \max_\alpha \bigg | K_\alpha^{(j)} - \sum_{\beta=1}^{\frac{n(n-1)}{2}} \omega_{\alpha\beta} K_\beta^{(x,t)} \bigg | \leq C \, 2^{-\frac{j}{2}}\] 
at the point $(x,t)$. 

For each point $(x,t) \in \partial \Omega^{(j)}$, we have $\max_\alpha |\langle K_\alpha^{(x,t)},\nu \rangle| \leq 2^{-j}$ at the point $(x,t)$. Consequently, $\max_\alpha |\langle K_\alpha^{(j)},\nu \rangle| \leq C \, 2^{-\frac{j}{2}}$ for all points $(x,t) \in \partial \Omega^{(j)}$ satisfying $\bar{t}-2^{\frac{j}{100}} \leq t \leq \bar{t}$. Finally, we note that, for each $t \in [\bar{t}-2^{\frac{j}{100}},\bar{t}]$, the point $p_t$ has distance at most $C \, 2^{\frac{j}{100}}$ from the point $p_{\bar{t}}$. This implies $\max_\alpha |\langle K_\alpha^{(j)},\nu \rangle| \leq C \, 2^{\frac{j}{100}}$ for all points $(x,t) \in \Omega^{(j)}$ with $t \in [\bar{t}-2^{\frac{j}{100}},\bar{t}]$. 

For each $\alpha \in \{1,\hdots,\frac{n(n-1)}{2}\}$, we define a function $f_\alpha^{(j)}: \Omega^{(j)} \to \mathbb{R}$ by 
\[f_\alpha^{(j)} := \exp(-2^{-\frac{j}{200}} (\bar{t}-t)) \, \frac{\langle K_\alpha^{(j)},\nu \rangle}{H - 2^{-\frac{j}{400}}}.\] 
Since $\max_\alpha |\langle K_\alpha^{(j)},\nu \rangle| \leq C \, 2^{-\frac{j}{2}}$ for all points $(x,t) \in \partial \Omega^{(j)}$ satisfying $\bar{t}-2^{\frac{j}{100}} \leq t \leq \bar{t}$, we conclude that 
\[\max_\alpha |f_\alpha^{(j)}(x,t)|\leq C\,2^{-\frac{j}{4}}\] 
for all points $(x,t) \in \partial \Omega^{(j)}$ satisfying $\bar{t}-2^{\frac{j}{100}} \leq t \leq \bar{t}$. The same estimate holds for all points $(x,t) \in \Omega^{(j)}$ with $t=\bar{t}-2^{\frac{j}{100}}$. On the other hand, the function $f_\alpha^{(j)}$ satisfies the evolution equation 
\[\frac{\partial}{\partial t} f_\alpha^{(j)} = \Delta f_\alpha^{(j)} + \frac{2}{H-2^{-\frac{j}{400}}} \, \langle \nabla H,\nabla f_\alpha^{(j)} \rangle - 2^{-\frac{j}{400}} \, \Big ( \frac{|A|^2}{H-2^{-\frac{j}{400}}} - 2^{-\frac{j}{400}} \Big ) \, f_\alpha^{(j)}.\] 
On the set $\Omega^{(j)}$, we have 
\[\frac{|A|^2}{H-2^{-\frac{j}{400}}} - 2^{-\frac{j}{400}} \geq \frac{1}{n} \, \frac{H^2}{H-2^{-\frac{j}{400}}} - 2^{-\frac{j}{400}} \geq \frac{1}{n} \, H - 2^{-\frac{j}{400}} \geq 0.\] 
By the maximum principle, we obtain 
\begin{align*} 
&\sup_{(x,t) \in \Omega^{(j)}, \, \bar{t}-2^{\frac{j}{100}} \leq t \leq \bar{t}} |f_\alpha^{(j)}(x,t)| \\ 
&\leq \max \bigg \{ \sup_{(x,t) \in \partial \Omega^{(j)}, \, \bar{t}-2^{\frac{j}{100}} \leq t \leq \bar{t}} |f_\alpha^{(j)}(x,t)|,\sup_{(x,t) \in \Omega^{(j)}, \, t=\bar{t}-2^{\frac{j}{100}}} |f_\alpha^{(j)}(x,t)| \bigg \} \\ 
&\leq C \, 2^{-\frac{j}{4}} 
\end{align*} 
for each $\alpha \in \{1,\hdots,\frac{n(n-1)}{2}\}$. Thus, we conclude that $\max_\alpha |\langle K_\alpha^{(j)},\nu \rangle| \leq C \, 2^{-\frac{j}{4}}$ for all points $(x,t) \in \Omega^{(j)}$ with $t=\bar{t}$. In particular, the distance of the point $p_{\bar{t}}$ from the axis of rotation of $\mathcal{K}^{(j)} = \{K_\alpha^{(j)}: 1 \leq \alpha \leq \frac{n(n-1)}{2}\}$ is bounded from above by a uniform constant which is independent of $j$. Sending $j \to \infty$, the sequence $\mathcal{K}^{(j)} = \{K_\alpha^{(j)}: 1 \leq \alpha \leq \frac{n(n-1)}{2}\}$ converges to a normalized set of rotation vector fields which are tangential along $M_{\bar{t}}$. This completes the proof of Theorem \ref{symmetry}. \\

Once we know that $M_t$ is rotationally symmetric for $-t$ sufficiently large, it follows from standard arguments that $M_t$ is rotationally symmetric for all $t \in (-\infty,0]$ (see \cite{Brendle-Choi}, Proposition 5.5).

\section{Uniqueness of ancient solutions with rotational symmetry}

\label{analysis.in.rotationally.symmetric.case}

Let $M_t$ be an ancient solution to the mean curvature flow in $\mathbb{R}^{n+1}$ which is strictly convex, uniformly two-convex, and noncollapsed. We may assume that $M_t$ is symmetric with respect to the $x_{n+1}$-axis. Let us write $M_t$ as a graph of a rotationally symmetric function $f$ on $\mathbb{R}^n$. The function $f$ satisfies the equation 
\[f_t= \frac{f_{rr}}{1+f_r^2}+\frac{n-1}{r}f_r.\]
Note that $f(r,t)$ may not be defined for all $r$.

Conversely, we may write the radius $r$ as a function of $(z,t)$, so that
\[f \big(r(z,t),t \big )=z.\]
Then $r(z,t)$ satisfies the equation 
\[r_t =\frac{r_{zz}}{1+r_z^2} -\frac{n-1}{r}.\]
Since $M_t$ is convex, we have 
\begin{align*} 
&r_z>0, &&r_t < 0 ,&& r_{zz} <0. 
\end{align*}
Without loss of generality, we assume that the tip of $M_0$ is at the origin. In other words, $f(0,0)=0$ and $r(0,0)=0$.

As in \cite{Brendle-Choi}, let $q_t = (0,\hdots,0,f(0,t))$ denote the tip of $M_t$, and let $H_{\text{\rm tip}}(t)$ denote the mean curvature of $M_t$ at the tip $q_t$. By the Harnack inequality \cite{Hamilton}, the function $t \mapsto H_{\text{\rm tip}}(t)$ is monotone increasing. Hence, the limit $\mathcal{H} := \lim_{t \to -\infty} H_{\text{\rm tip}}(t)$ exists. Using results in Section \ref{barrier}, we obtain $|q_t| \geq c \, (-t)$ for $-t$ sufficiently large. This gives $\mathcal{H} > 0$.

We first prove that $f_t(r,t)$ is monotone increasing in $t$.

\begin{proposition}
\label{consequence.of.harnack}
We have $f_{tt}(r,t)\geq 0$ everywhere.
\end{proposition}

\textbf{Proof.} 
This is a consequence of Hamilton's Harnack inequality for mean curvature flow \cite{Hamilton}. See \cite{Brendle-Choi} for details. \\

We next show that $f_t(r,t)$ is bounded from below.

\begin{proposition}
\label{lower.bound.for.f_t}
We have $f_t(r,t) \geq \mathcal{H}$ at each point in space-time. Moreover, for each $r_0 > 0$,
\[\lim_{t \to -\infty} \sup_{r \leq r_0} f_t(r,t) = \mathcal{H}.\]
\end{proposition}

\textbf{Proof.}
Consider an arbitrary sequence of times $t_j \to -\infty$, and define $M_t^{(j)} := M_{t+t_j}-q_{t_j}$. Since the flows $M_t^{(j)}$ have uniformly bounded curvature, the sequence $M_t^{(j)}$ converges in $C_{\text{\rm loc}}^\infty$ to a smooth eternal solution, which is rotationally symmetric. At each point in time, the mean curvature at the tip of the limit solution equals $\mathcal{H}$. Consequently, we are in the equality case in the Harnack inequality. By \cite{Hamilton}, the limit solution must be a self-similar translator which is moving with speed $\mathcal{H}$. This gives 
\[\lim_{j \to \infty} \sup_{r \leq r_0} |f_t(r,t_j) - \mathcal{H}| = 0\]
for every $r_0>0$. Since $f_{tt}(r,t) \geq 0$ by Proposition \ref{consequence.of.harnack}, it follows that $f_t(r,t) \geq \mathcal{H}$ for all $r$ and $t$. \\

In the next step, we show that $f_t(r,t)$ is monotone increasing in $r$.

\begin{proposition}
\label{monotonicity.in.r}
We have $f_{tr}(r,t) \geq 0 $ everywhere. 
\end{proposition}

\textbf{Proof.} 
Consider a time $t_0$ and a radius $r_0$ such that $f(r_0,t_0)$ is defined. Moreover, let $t_j$ be a sequence of times such that $t_j \to -\infty$. For $j$ large, we denote by $Q_j$ the parabolic cylinder $Q_j=\{x_1^2+\hdots+x_n^2 \leq r_0^2, t \in [t_j,t_0]\}$. It follows from the evolution equations for $H$ and $\langle \omega,\nu \rangle$ that the maximum $\sup_{Q_j} H \, \langle \omega,\nu \rangle^{-1}$ must be attained on the parabolic boundary of $Q_j$. This gives 
\begin{align*} 
&\sup_{x_1^2+\hdots+x_n^2 \leq r_0^2, t=t_0} H \, \langle \omega,\nu \rangle^{-1} \\ 
&\leq \max \Big \{\sup_{x_1^2+\hdots+x_n^2 = r_0^2, t_j \leq t \leq t_0} H \, \langle \omega,\nu \rangle^{-1}, \, \sup_{x_1^2+\hdots+x_n^2 \leq r_0^2, t=t_j} H \, \langle \omega,\nu \rangle^{-1} \Big \}. 
\end{align*}
Since $f_t(r,t) = H \, \langle \omega,\nu \rangle^{-1}$, it follows that 
\begin{align*} 
\sup_{r \leq r_0} f_t(r,t_0) 
&\leq \max \Big \{ \sup_{t_j \leq t \leq t_0} f_t(r_0,t), \, \sup_{r \leq r_0} f_t(r,t_j) \Big \} \\ 
&= \max \Big \{ f_t(r_0,t_0), \, \sup_{r \leq r_0} f_t(r,t_j) \Big \}. 
\end{align*}
Note that in the last step we have used Proposition \ref{consequence.of.harnack}. Finally, we pass to the limit as $j \to \infty$. Using Proposition \ref{lower.bound.for.f_t}, we obtain $\sup_{r \leq r_0} f_t(r,t_0) \leq \max \{f_t(r_0,t_0),\mathcal{H}\} = f_t(r_0,t_0)$. This completes the proof of Proposition \ref{monotonicity.in.r}. \\

By assumption, $M_t$ is strictly convex, uniformly two-convex, and noncollapsed. Moreover, $H_{\text{\rm tip}}(t)$ is bounded from below by $\mathcal{H}$. Hence, there exists a small constant $\varepsilon_0 \in (0,\frac{1}{20})$ and a decreasing function $\Lambda: (0,\varepsilon_0] \to \mathbb{R}$ such that given any $\varepsilon \in (0,\varepsilon_0]$, if $|\bar{x}-q_t| \geq \Lambda(\varepsilon)$, then $(\bar{x}, \bar{t})$ is a center of $\varepsilon$-neck (cf. \cite{Haslhofer-Kleiner2}, Proposition 3.1). We recall three estimates from \cite{Brendle-Choi}. These results were stated for $n=2$ in \cite{Brendle-Choi}, but the arguments carry over directly to higher dimensions.

\begin{lemma}
\label{bound.for.rr_z}
On every $\varepsilon_0$-neck, $rr_z=\frac{r}{f_r} \leq (1+2\varepsilon_0)(n-1)\mathcal{H}^{-1}$. 
\end{lemma}

\textbf{Proof.} 
This follows from the inequality $f_t \geq \mathcal{H}$ established in Proposition \ref{lower.bound.for.f_t}. See \cite{Brendle-Choi}, Lemma 6.4, for details. \\

\begin{lemma}
\label{interior.estimate.for.r}
There exists a constant $C_0 \geq 1$ such that $r^m \, \big | \frac{\partial^m}{\partial z^m} r \big | \leq C_0$ holds for $m=1,2,3$ at the center of $\varepsilon_0$-necks with $r \geq 1$. 
\end{lemma}

\textbf{Proof.} 
Analogous to Lemma 6.5 in \cite{Brendle-Choi}. \\

\begin{proposition}
\label{estimate.for.r_zz}
There exist constants $C_1$ and $C_2$ with the following property. If $r \geq C_1$, then $0 \leq -r_{zz}(z,t) \leq C_2 r(z,t)^{-\frac{5}{2}}$.
\end{proposition}

\textbf{Proof.} 
Analogous to Proposition 6.6 in \cite{Brendle-Choi}. \\

For each $z<0$, we define a real number $\mathcal{T}(z)$ by
\begin{align*}
&r(z,t)>0 \quad \text{for}\quad t < \mathcal{T}(z), && \lim_{t\to \mathcal{T}(z)} r(z,t)=0.
\end{align*}
In other words, if $t=\mathcal{T}(z)$, then $f(0,t) = z$, and the tip of $M_t$ is located at $(0,\hdots,0,z)$. The following result allows us to estimate $r(z,t)$ in terms of $\mathcal{T}(z)-t$.

\begin{corollary}
\label{cor}
We have 
\[2(n-1) \, [\mathcal{T}(z)-t] \leq r(z,t)^2 \leq 2(n-1) \, [\mathcal{T}(z)-t] + 8C_2[\mathcal{T}(z)-t]^{\frac{1}{4}}+C_1^2\] 
if $z<0$ and $r(z,t)$ is sufficiently large.
\end{corollary}

\textbf{Proof.} Let us fix a point $(\bar{z},\bar{t})$. The inequality $(r^2+2(n-1)t)_t = \frac{2rr_{zz}}{1+r_z^2} < 0$ implies 
\[r(\bar{z},\bar{t})^2 \geq 2(n-1) \, [\mathcal{T}(\bar{z})-\bar{t}].\] 
Moreover, if $r \geq C_1$, then $(r^2+2(n-1)t)_t = \frac{2rr_{zz}}{1+r_z^2} \geq -2C_2 r^{-\frac{3}{2}}$ by Proposition \ref{estimate.for.r_zz}. Let $\tilde{t} \leq \mathcal{T}(\bar{z})$ be chosen so that $r(\bar{z},\tilde{t})=C_1$. Then
\begin{align*}
r(\bar{z},\bar{t})^2
&= C_1^2+2(n-1)(\tilde{t}-\bar{t}) - \int_{\bar{t}}^{\tilde{t}} (r(\bar{z},t)^2+2(n-1)t)_t \, dt \\
&\leq C_1^2+2(n-1)(\tilde{t}-\bar{t})+2C_2\int_{\bar{t}}^{\tilde{t}} r(\bar{z},t)^{-\frac{3}{2}} \, dt \\ 
&\leq C_1^2+2(n-1)(\tilde{t}-\bar{t})+2C_2\int_{\bar{t}}^{\tilde{t}} [\mathcal{T}(\bar{z})-t]^{-\frac{3}{4}} \, dt \\
&\leq C_1^2+2(n-1)(\tilde{t}-\bar{t})-8C_2[\mathcal{T}(\bar{z})-\tilde{t}]^{\frac{1}{4}}+8C_2[\mathcal{T}(\bar{z})-\bar{t}]^{\frac{1}{4}} \\ 
&\leq C_1^2 + 2(n-1) \, [\mathcal{T}(\bar{z})-\bar{t}] + 8C_2[\mathcal{T}(\bar{z})-\bar{t}]^{\frac{1}{4}}.
\end{align*} 
This proves the assertion. \\

\begin{lemma}
\label{lower.bound.for.rr_z.at.z=0}
Let $\delta$ be an arbitrary positive real number. Then 
\[r(0,t)r_z(0,t) \geq (n-1) \, (\mathcal{H}^{-1}-\delta)\]
provided that $-t$ is sufficiently large.
\end{lemma}

\textbf{Proof.} 
We may assume that $R=r(0,t) \geq C_1$. Then every point $(x,t)$ with $x_{n+1}=0$ lies at the center of an $\varepsilon_0$-neck. Consequently, $|r(z,t)-R| \leq \varepsilon_0 R$ for $|z| \leq 2R$. Moreover, we have $rr_z \leq (1+2\varepsilon_0) \, (n-1) \, \mathcal{H}^{-1}$ by Lemma \ref{bound.for.rr_z}, and $|(rr_{z})_z| = |rr_{zz}+r_z^2| \leq C_3R^{-\frac{3}{2}}$ by Proposition \ref{estimate.for.r_zz}. Hence, if $R^{\frac{1}{2}}\geq 4C_3\delta^{-1}$, then we have
\[|r(z,t)r_z(z,t)-r(0,t)r_z(0,t)| \leq 2C_3R^{-\frac{1}{2}} \leq \frac{\delta}{2}\]
for all $z \in [-2R,2R]$.

Using Corollary \ref{cor}, we obtain 
\[r(-R,t)^2 \geq 2(n-1) \, [\mathcal{T}(-R)-t],\] 
\[r(-2R,t)^2 \geq 2(n-1) \, [\mathcal{T}(-2R)-t],\] 
and 
\begin{align*} 
r(-2R,t)^2 
&\leq 2(n-1) \, [\mathcal{T}(-2R)-t] + 8C_2 [\mathcal{T}(-2R)-t]^{\frac{1}{4}}+C_1^2 \\ 
&\leq 2(n-1) \, [\mathcal{T}(-2R)-t] + 8C_2 r(-2R,t)^{\frac{1}{2}}+C_1^2 \\  
&\leq 2(n-1) \, [\mathcal{T}(-2R)-t] + 8C_2 R^{\frac{1}{2}}+C_1^2,
\end{align*}
where in the last step we have used the inequality $r(-2R,t) \leq r(0,t) = R$. From this, we deduce that  
\[r(-R,t)^2-r(-2R,t)^2 \geq 2(n-1) \, [\mathcal{T}(-R)-\mathcal{T}(-2R)]-8C_2R^{\frac{1}{2}}-C_1^2.\] 
On the other hand, using the identity $f_t(0,t) = H_{\text{\rm tip}}(t)$, we obtain $\frac{d}{dz} \mathcal{T}(z) = H_{\text{\rm tip}}(\mathcal{T}(z))^{-1}$. Hence, if $R$ is sufficiently large, then $\frac{d}{dz} \mathcal{T}(z) \geq \mathcal{H}^{-1} - \frac{\delta}{4}$ for $z \in [-2R,-R]$. This gives 
\[\mathcal{T}(-R) - \mathcal{T}(-2R) \geq \Big ( \mathcal{H}^{-1}-\frac{\delta}{4} \Big ) \, R\] 
if $R$ is sufficiently large. Putting these facts together, we obtain 
\[r(-R,t)^2-r(-2R,t)^2 \geq 2(n-1) \, \Big ( \mathcal{H}^{-1}-\frac{\delta}{2} \Big ) \, R,\] 
hence 
\[\sup_{z \in [-2R,-R]} r(z,t) r_z(z,t) \geq (n-1) \, \Big ( \mathcal{H}^{-1}-\frac{\delta}{2} \Big )\] 
if $-t$ is sufficiently large. Thus, 
\[r(0,t)r_z(0,t) \geq (n-1) \, (\mathcal{H}^{-1}-\delta)\] 
if $-t$ is sufficiently large. \\

\begin{proposition}
\label{liminf}
Given $\delta>0$, there exists a time $\bar{t} \in (-\infty,0]$ (depending on $\delta$) such that 
\begin{align*}
r(z,t)r_z(z,t) \geq (n-1) \, (\mathcal{H}^{-1}-2\delta),
\end{align*}
holds for all $z \geq 0$ and $t \leq \bar{t}$.
\end{proposition}

\textbf{Proof.} 
Let $\psi(z,t)$ denote the solution of the Dirichlet problem for the one-dimensional heat equation on the half line with initial condition $\lim_{t \to 0} \psi(z,t) = 1$ (see \cite{Brendle-Choi}, Proposition 6.9). By Lemma \ref{lower.bound.for.rr_z.at.z=0}, we can find a time $\bar{t}$ so that $r(0,t)r_z(0,t) \geq (n-1) \, (\mathcal{H}^{-1}-\delta)$ for $t \leq \bar{t}$. Moreover, we can arrange that $r(z,t) \geq C_1+C_2$ for $z \geq 0$ and $t \leq \bar{t}$. Given any $s < \bar{t}$, we define a function $\psi^{\delta,s}(z,t)$ by 
\[\psi^{\delta,s}(z,t) = (n-1) \, (\mathcal{H}^{-1}-2\delta -\mathcal{H}^{-1} \, \psi(2z,t-s))\] 
for $t \in (s,\bar{t}]$. We will show that $rr_z > \psi^{\delta,s}$ for all $z \geq 0$ and all $t \in (s,\bar{t}]$.  

It is straightforward to verify that $r(0,t) r_z(0,t) \geq (n-1) \, (\mathcal{H}^{-1}-\delta) > \limsup_{z \to 0} \psi^{\delta,s}(z,t)$ for each $t \in (s,\bar{t}]$; $\liminf_{z \to \infty} r(z,t) r_z(z,t) \geq 0 > \limsup_{z \to \infty} \psi^{\delta,s} (z,t)$ for each $t \in (s,\bar{t}]$; and $r(z,s) r_z(z,s) \geq 0 > \limsup_{t \to s} \psi^{\delta,s} (z,t)$ for each $z>0$. 

On the other hand, for $z > 0$ and $t \in (s,\bar{t}]$, we have , we have $1+rr_{zz} \geq 0$, hence 
\[(rr_z)_t = \frac{(rr_z)_{zz}}{1+r_z^2} - \frac{2r_zr_{zz} (1+r_z^2+rr_{zz})}{(1+r_z^2)^2} > \frac{(rr_z)_{zz}}{1+r_z^2}.\] 
Moreover, 
\[(\psi^{\delta,s})_t = \frac{1}{4} \, (\psi^{\delta,s})_{zz} \leq \frac{(\psi^{\delta,s})_{zz}}{1+r_z^2}.\] 
Using the maximum principle, we conclude that $rr_z > \psi^{\delta,s}$ for all $z \geq 0$ and all $t \in (s,\bar{t}]$. Sending $s \to -\infty$ gives the desired result. \\

\begin{corollary}
\label{f.defined.everywhere}
We can find a time $T \in (-\infty,0]$ such that $r(z,t)^2 \geq (n-1) \, \mathcal{H}^{-1} \, z$ for all $z \geq 0$ and $t \leq T$. In particular, if $t \leq T$, then the function $f(r,t)$ is defined for all $r \in [0,\infty)$.
\end{corollary}

\textbf{Proof.} By Proposition \ref{liminf}, we can find a time $T \in (-\infty,0]$ such that $r(z,t)r_z(z,t) \geq \frac{1}{2} \, (n-1) \, \mathcal{H}^{-1}$ for all $z \geq 0$ and all $t \leq T$. If we integrate over $z$, the assertion follows. \\

\begin{proposition}
\label{limit.of.rr_z}
For each $t \leq T$, we have $\lim_{z \to \infty} r(z,t)r_z(z,t) = (n-1) \, \mathcal{H}^{-1}$.
\end{proposition}

\textbf{Proof.}
It follows from Lemma \ref{bound.for.rr_z} that  $\limsup_{z \to \infty} r(z,t)r_z(z,t) \leq (n-1) \, \mathcal{H}^{-1}$. Hence, it suffices to prove that $\liminf_{z \to \infty} r(z,t) r_z(z,t) \geq (n-1) \, \mathcal{H}^{-1}$ for each $t \leq T$. By Proposition \ref{liminf}, we know that $\liminf_{z \to \infty} r(z,t) r_z(z,t) \geq (n-1) \, (\mathcal{H}^{-1}-2\delta)$ if $-t$ is sufficiently large. On the other hand, Lemma \ref{interior.estimate.for.r} gives 
\[|(rr_z)_t|=|r_tr_z+rr_{zt}| = \Big | \frac{r_zr_{zz}}{1+r_z^2}+\frac{rr_{zzz}}{1+r_z^2}-\frac{2rr_zr_{zz}^2}{(1+r_z^2)^2} \Big | \leq \frac{4C_0^3}{r^2}\] 
for $r \geq C_1$. Using Corollary \ref{f.defined.everywhere}, we conclude that the quantity $\liminf_{z \to \infty} r(z,t)r_z(z,t)$ has the same value for each $t \leq T$. Putting these facts together, we conclude that $\liminf_{z \to \infty} r(z,t) r_z(z,t) \geq (n-1) \, \mathcal{H}^{-1}$ for each $t \leq T$. \\

\begin{theorem}
\label{soliton}
For each $t \leq T$, the solution $M_t$ is a rotationally symmetric translating soliton.
\end{theorem}

\textbf{Proof.} 
Since $rr_z =\frac{r}{f_r}$, Proposition \ref{limit.of.rr_z} implies
\[\lim_{r\to \infty} \frac{f_r(r,t)}{r}=\frac{1}{n-1} \, \mathcal{H}\]
for each $t \leq T$. Moreover, Lemma \ref{interior.estimate.for.r} gives $\limsup_{z \to \infty} r^2 \, (-r_{zz}) < \infty$, hence $\limsup_{r \to \infty} \frac{r^2 \, f_{rr}}{f_r^3} < \infty$. Consequently, 
\[\lim_{r \to \infty} \frac{f_{rr}(r,t)}{1+f_r(r,t)^2} = 0.\] 
Using the evolution equation for $f(r,t)$, we obtain 
\[\lim_{r\to \infty} f_t(r,t) = \lim_{r \to \infty} \frac{(n-1) \, f_r(r,t)}{r}=\mathcal{H}\] 
for each $t \leq T$. Using Proposition \ref{monotonicity.in.r}, we conclude that $f_t(r,t) \leq \mathcal{H}$ for all $r \geq 0$ and all $t \leq T$. Therefore, Proposition \ref{lower.bound.for.f_t} gives $f_t(r,t)=\mathcal{H}$ for all $r \geq 0$ and all $t \leq T$. Consequently, $M_t$ is a translating soliton for each $t \leq T$. \\

Once we know that $M_t$ is a translating soliton for $-t$ sufficiently large, it follows from standard arguments that $M_t$ is a translating soliton for all $t \in (-\infty,0]$ (see \cite{Brendle-Choi}, Proposition 6.14). This completes the proof of Theorem \ref{main.thm}.

\end{document}